%
%
%
%
%
\newif\ifsect\newif\iffinal
\secttrue\finaltrue
\def\strutdepth{\dp\strutbox}

\def\lsimb#1{\vadjust{\vtop to0pt{\baselineskip\strutdepth\vss
	\llap{\ttt\string #1\ }\null}}}
\def\rsimb#1{\vadjust{\vtop to0pt{\baselineskip\strutdepth\vss
	\line{\kern\hsize\rlap{\ttt\ \string #1}}\null}}}
\def\ssect #1. {\bigbreak\indent{\bf #1.}\enspace\message{#1}}
\def\smallsect #1. #2\par{\bigbreak\noindent{\bf #1.}\enspace{\bf #2}\par
	\global\parano=#1\global\eqnumbo=1\global\thmno=1\global\esno=0
	\global\remno=0\global\remno=0\global\defno=0
	\nobreak\smallskip\nobreak\noindent\message{#2}}
\def\thm #1: #2{\medbreak\noindent{\bf #1:}\if(#2\thmp\else\thmn#2\fi}
\def\thmp #1) { (#1)\thmn{}}
\def\thmn#1#2\par{\enspace{\sl #1#2}\par
        \ifdim\lastskip<\medskipamount \removelastskip\penalty 55\medskip\fi}
\def\qedn{\thinspace\null\nobreak\hfill\hbox{\vbox{\kern-.2pt\hrule height.2pt 
depth.2pt\kern-.2pt\kern-.2pt \hbox to2.5mm{\kern-.2pt\vrule width.4pt
\kern-.2pt\raise2.5mm\vbox to.2pt{}\lower0pt\vtop to.2pt{}\hfil\kern-.2pt
\vrule width.4pt\kern-.2pt}\kern-.2pt\kern-.2pt\hrule height.2pt depth.2pt
\kern-.2pt}}\par\medbreak}
\def\pf{\ifdim\lastskip<\smallskipamount \removelastskip\smallskip\fi
        \noindent{\sl Proof\/}:\enspace}
\def\itm#1{\par\indent\llap{\rm #1\enspace}\ignorespaces}

\def\bar#1{\overline{#1}}
\def\forclose#1{\hfil\llap{$#1$}\hfilneg}
\def\newforclose#1{
	\ifsect\xdef #1{(\number\parano.\number\eqnumbo)}\else
	\xdef #1{(\number\eqnumbo)}\fi
	\hfil\llap{$#1$}\hfilneg
	\global \advance \eqnumbo by 1
	\iffinal\else\rsimb#1\fi}
\def\newforclosea#1{
	\ifsect\xdef #1{{\rm(\the\Apptok.\number\eqnumbo)}}\else
	\xdef #1{(\number\eqnumbo)}\fi
	\hfil\llap{$#1$}\hfilneg
	\global \advance \eqnumbo by 1
	\iffinal\else\rsimb#1\fi}
\def\forevery#1#2$${\displaylines{\let\eqno=\forclose\let\neweqa=\newforclosea
        \let\neweq=\newforclose\hfilneg\rlap{$\qqquad\forall#1$}\hfil#2\cr}$$}
\newcount\parano
\newcount\eqnumbo
\newcount\thmno
\newcount\versiono
\versiono=0
\def\neweqt#1$${\xdef #1{(\number\parano.\number\eqnumbo)}
	\eqno #1$$
	\iffinal\else\rsimb#1\fi
	\global \advance \eqnumbo by 1}
\def\newthmt#1 #2: #3{\xdef #2{\number\parano.\number\thmno}
	\global \advance \thmno by 1
	\medbreak\noindent
	\iffinal\else\lsimb#2\fi
	{\bf #1 #2:}\if(#3\thmp\else\thmn#3\fi}
\def\neweqf#1$${\xdef #1{(\number\eqnumbo)}
	\eqno #1$$
	\iffinal\else\rlap{$\smash{\hbox{\hfilneg\string#1\hfilneg}}$}\fi
	\global \advance \eqnumbo by 1}
\def\newthmf#1 #2: #3{\xdef #2{\number\thmno}
	\global \advance \thmno by 1
	\medbreak\noindent
	\iffinal\else\llap{$\smash{\hbox{\hfilneg\string#1\hfilneg}}$}\fi
	{\bf #1 #2:}\if(#3\thmp\else\thmn#3\fi}
\def\inizia{\ifsect\let\neweq=\neweqt\else\let\neweq=\neweqf\fi
\ifsect\let\newthm=\newthmt\else\let\newthm=\newthmf\fi}
\def\bititolo{\empty}
\gdef\begin #1 #2\par{\xdef\titolo{#2}
\ifsect\let\neweq=\neweqt\else\let\neweq=\neweqf\fi
\ifsect\let\newthm=\newthmt\else\let\newthm=\newthmf\fi
\centerline{\titlefont\titolo}
\if\bititolo\empty\else\medskip\centerline{\titlefont\bititolo}\fi
\bigskip
\centerline{\bigfont \autore}
\if\istituto!\else\bigskip
\centerline{\istituto}
\centerline{\indirizzo}
\centerline{\email}\fi
\medskip
\centerline{#1~\anno}
\bigskip\bigskip
\ifsect\else\global\thmno=1\global\eqnumbo=1\fi}
\font\titlefont=cmssbx10 scaled \magstep1
\font\bigfont=cmr12

\font\ttt=cmtt10 at 10truept
\font\eightrm=cmr8

\let\sc=\smallcaps

\font\bbr=msbm10
\font\sbbr=msbm7 
\font\ssbbr=msbm5
\def\ca #1{{\cal #1}}
\nopagenumbers
\binoppenalty=10000
\relpenalty=10000
\newfam\amsfam
\textfont\amsfam=\bbr \scriptfont\amsfam=\sbbr \scriptscriptfont\amsfam=\ssbbr
\let\de=\partial
\def\eps{\varepsilon}
\def\phe{\varphi}

\def\Hom{\mathop{\rm Hom}\nolimits}
\def\End{\mathop{\rm End}\nolimits}

\def\Re{\mathop{\rm Re}\nolimits}
\def\Im{\mathop{\rm Im}\nolimits}

\def\rk{\mathop{\rm rk}\nolimits}
\def\id{\mathop{\rm id}\nolimits}
\mathchardef\void="083F
\def\Z{{\mathchoice{\hbox{\bbr Z}}{\hbox{\bbr Z}}{\hbox{\sbbr Z}}
{\hbox{\sbbr Z}}}}
\def\R{{\mathchoice{\hbox{\bbr R}}{\hbox{\bbr R}}{\hbox{\sbbr R}}
{\hbox{\sbbr R}}}}
\def\C{{\mathchoice{\hbox{\bbr C}}{\hbox{\bbr C}}{\hbox{\sbbr C}}
{\hbox{\sbbr C}}}}
\def\N{{\mathchoice{\hbox{\bbr N}}{\hbox{\bbr N}}{\hbox{\sbbr N}}
{\hbox{\sbbr N}}}}
\def\P{{\mathchoice{\hbox{\bbr P}}{\hbox{\bbr P}}{\hbox{\sbbr P}}
{\hbox{\sbbr P}}}}
\def\Q{{\mathchoice{\hbox{\bbr Q}}{\hbox{\bbr Q}}{\hbox{\sbbr Q}}
{\hbox{\sbbr Q}}}}

\def\qqquad{\quad\qquad}

\newcount\notitle
\notitle=1
\headline={\ifodd\pageno\rhead\else\lhead\fi}
\def\rhead{\ifnum\pageno=\notitle\iffinal\hfill\else\hfill\tt Version 
\the\versiono; \the\day/\the\month/\the\year\fi\else\hfill\eightrm\titolo\hfill
\folio\fi}
\def\lhead{\ifnum\pageno=\notitle\hfill\else\eightrm\folio\hfill\autore\hfill
\fi}
\def\autore{Marco Abate}
\output={\plainoutput}
\newbox\bibliobox
\def\setref #1{\setbox\bibliobox=\hbox{[#1]\enspace}
	\parindent=\wd\bibliobox}
\def\biblap#1{\noindent\hang\rlap{[#1]\enspace}\indent\ignorespaces}
\def\art#1 #2: #3! #4! #5 #6 #7-#8 \par{\biblap{#1}#2: {\sl #3\/}.
	#4 {\bf #5} (#6)\if.#7\else, \hbox{#7--#8}\fi.\par\smallskip}
\def\book#1 #2: #3! #4 \par{\biblap{#1}#2: {\bf #3.} #4.\par\smallskip}
\def\coll#1 #2: #3! #4! #5 \par{\biblap{#1}#2: {\sl #3\/}. In {\bf #4,} 
#5.\par\smallskip}
\def\pre#1 #2: #3! #4! #5 \par{\biblap{#1}#2: {\sl #3\/}. #4, #5.\par\smallskip}
%
%
\let\newthm=\newthmt
\let\neweq=\neweqt
\newcount\esno\newcount\defno\newcount\remno
\def\Def #1\par{\global \advance \defno by 1
    \medbreak
{\bf Definition \the\parano.\the\defno:}\enspace #1\par
\ifdim\lastskip<\medskipamount \removelastskip\penalty 55\medskip\fi}
\def\Rem #1\par{\global \advance \remno by 1
    \medbreak
{\bf Remark \the\parano.\the\remno:}\enspace #1\par
\ifdim\lastskip<\medskipamount \removelastskip\penalty 55\medskip\fi}
\def\Es #1\par{\global \advance \esno by 1
    \medbreak
{\sc Example \the\parano.\the\esno:}\enspace #1\par
\ifdim\lastskip<\medskipamount \removelastskip\penalty 55\medskip\fi}
\def\istituto{Dipartimento di Matematica, Universit\`a di Pisa}
\def\indirizzo{Largo Pontecorvo 5, 56127 Pisa, Italy}
\def\anno{2008}
\def\email{E-mail: abate@dm.unipi.it}
%
%
%
%
\begin {November} Discrete holomorphic local dynamical systems 

\smallsect 1. Introduction

Let us begin by defining the main object of study in this survey.

\Def Let $M$ be a complex manifold, and $p\in M$. A {\sl
(discrete) holomorphic local dynamical system} at~$p$ is a holomorphic map
$f\colon U\to M$ such that~$f(p)=p$, where $U\subseteq M$ is an open
neighbourhood of~$p$; we shall also assume that $f\not\equiv\id_U$. We shall
denote by~$\End(M,p)$ the set of holomorphic local dynamical systems at~$p$. 

\Rem Since we are mainly concerned with the behavior of~$f$ nearby~$p$, we shall
sometimes replace~$f$ by its restriction to some suitable open
neighbourhood of~$p$. It is possible to formalize this fact by using germs of
maps and germs of sets at~$p$, but for our purposes it will be enough
to use a somewhat less formal approach.

\Rem In this survey we shall never have the occasion of discussing continuous
holomorphic dynamical systems (i.e., holomorphic foliations). So from now on all
dynamical systems in this paper will be discrete, except where explicitly noted
otherwise.

To talk about the dynamics of an $f\in\End(M,p)$ we need to define the iterates
of~$f$. If $f$ is defined on the set~$U$, then the second iterate $f^2=f\circ f$
is defined on~$U\cap f^{-1}(U)$ only, which still is an open neighbourhood
of~$p$. More generally, the $k$-th iterate $f^k=f\circ f^{k-1}$ is defined
on~$U\cap f^{-1}(U)\cap\cdots\cap f^{-(k-1)}(U)$. This 
suggests the next definition:

\Def Let $f\in\End(M,p)$ be a holomorphic local dynamical system defined
on an open set $U\subseteq M$. Then the {\sl stable set}~$K_f$ of~$f$ is
$$
K_f=\bigcap_{k=0}^\infty f^{-k}(U)\;.
$$
In other words, the stable set of~$f$ is
the set of all points~$z\in U$ such that the {\sl orbit} $\{f^k(z)\mid
k\in\N\}$ is well-defined. If $z\in U\setminus K_f$, we shall say that $z$ (or
its orbit) {\sl escapes} from~$U$. 

Clearly, $p\in K_f$, and so the stable set is never empty (but it can happen
that $K_f=\{p\}$; see the next section for an example). 
Thus the first natural question in local holomorphic dynamics is:
\smallskip
\item{(Q1)} {\it What is the topological structure of~$K_f$?}
\smallskip
\noindent For instance, when does $K_f$ have non-empty interior? As we shall see
in Proposition~4.1, holomorphic local dynamical systems such that $p$ belongs to the
interior of the stable set enjoy special properties.

\Rem Both the definition of stable set and Question~1 (as well as several other
definitions and questions we shall see later on) are topological in character;
we might state them for local dynamical systems which are continuous only. As
we shall see, however, the {\it answers} will strongly depend on the
holomorphicity of the dynamical system.

\Def Given $f\in\End(M,p)$, a set $K\subseteq M$ is {\sl completely 
$f$-invariant} if $f^{-1}(K)=K$ (this implies, in particular, that 
$K$ is {\sl $f$-invariant,} that is $f(K)\subseteq K$).

Clearly, the stable set~$K_f$ is completely $f$-invariant.
Therefore the pair $(K_f,f)$ is a discrete dynamical system in the usual
sense, and so the second natural question in local holomorphic dynamics is
\smallskip
\item{(Q2)} {\it What is the dynamical structure of~$(K_f,f)$?}
\smallskip
\noindent For instance, what is the asymptotic behavior of the orbits? Do
they converge to~$p$, or have they a chaotic behavior? Is there a dense
orbit? Do there exist proper $f$-invariant subsets, that is sets
$L\subset K_f$ such that~$f(L)\subseteq L$? If they do exist, what is the
dynamics on them?

To answer all these questions, the most efficient way is to replace $f$ by a
``dynamically equivalent" but simpler (e.g., linear) map~$g$. In our context,
``dynamically equivalent" means ``locally conjugated"; and we have at least
three kinds of conjugacy to consider.

\Def Let $f_1\colon U_1\to M_1$ and $f_2\colon U_2\to M_2$ be two holomorphic local
dynamical systems at~$p_1\in M_1$ and~$p_2\in M_2$ respectively. We shall say
that~$f_1$ and~$f_2$ are {\sl holomorphically} (respectively, {\sl
topologically\/}) {\sl locally conjugated} if there are open neighbourhoods
$W_1\subseteq U_1$ of~$p_1$, $W_2\subseteq U_2$ of~$p_2$, and a
biholomorphism (respectively, a homeomorphism) $\phe\colon W_1\to W_2$ with
$\phe(p_1)=p_2$ such that 
$$
f_1=\phe^{-1}\circ f_2\circ\phe
\qquad\hbox{on}\qquad\phe^{-1}\bigl(
W_2\cap f_2^{-1}(W_2)\bigr)=W_1\cap f_1^{-1}(W_1)\;.
$$  

If $f_1\colon U_1\to M_1$ and $f_2\colon U_2\to M_2$ are locally conjugated, 
in particular we have
$$
\forevery{k\in\N}\qquad\qquad f_1^k=\phe^{-1}\circ f_2^k\circ\phe
\quad\hbox{on}\quad\phe^{-1}\bigl(
W_2\cap\cdots\cap f_2^{-(k-1)}(W_2)\bigr)=W_1\cap\cdots\cap f_1^{-(k-1)}(W_1)\;,
$$
and thus 
$$
K_{f_2|_{W_2}}=\phe(K_{f_1|_{W_1}})\;.
$$
So the local dynamics of~$f_1$
about~$p_1$ is to all purposes equivalent to the local dynamics of~$f_2$
about~$p_2$.

\Rem Using local coordinates centered at~$p\in M$ it is easy to show that any
holomorphic local dynamical system at~$p$ is holomorphically locally conjugated
to a holomorphic local dynamical system at~$O\in\C^n$, where $n=\dim M$.

Whenever we have an equivalence relation in a class of objects, there are
classification problems. So the third natural question in local
holomorphic dynamics is
\smallskip
\item{(Q3)} {\it Find a (possibly small) class $\ca F$ of holomorphic local
dynamical systems at~$O\in\C^n$ such that every holomorphic local dynamical
system~$f$ at a point in an $n$-dimensional complex manifold is holomorphically
(respectively, topologically) locally conjugated to a (possibly) unique element
of~$\ca F$, called the {\sl holomorphic} (respectively, {\sl topological\/})
{\sl normal form} of~$f$.}
\smallskip
\noindent Unfortunately, the holomorphic classification is often too complicated
to be practical; the family~$\ca F$ of normal forms might be uncountable. A
possible replacement is looking for invariants instead of normal forms:
\smallskip
\item{(Q4)} {\it Find a way to associate a (possibly small) class of (possibly
computable) objects, called {\sl invariants,} to any holomorphic local
dynamical system~$f$ at~$O\in\C^n$ so that two holomorphic local dynamical
systems at~$O$ can be holomorphically conjugated only if they have the
same invariants. The class of invariants is furthermore said {\sl complete} if
two holomorphic local dynamical systems at~$O$ are holomorphically conjugated
if and only if they have the same invariants.}
\smallskip
\noindent As remarked before, up to now all the questions we asked made sense
for topological local dynamical systems; the next one instead makes sense only
for holomorphic local dynamical systems.

A holomorphic local dynamical system at~$O\in\C^n$ is clearly given by an
element of~$\C_0\{z_1,\ldots,z_n\}^n$, the space of $n$-uples of converging
power series in~$z_1,\ldots,z_n$ without constant terms. The
space~$\C_0\{z_1,\ldots,z_n\}^n$ is a subspace of the
space~$\C_0[[z_1,\ldots,z_n]]^n$ of $n$-uples of formal power series without
constant terms. An element $\Phi\in \C_0[[z_1,\ldots,z_n]]^n$ has an inverse
(with respect to composition) still belonging to~$\C_0[[z_1,\ldots,z_n]]^n$ if
and only if its linear part is a linear automorphism of~$\C^n$. 

\Def We shall 
say that two holomorphic local dynamical systems
$f_1$,~$f_2\in\C_0\{z_1,\ldots,z_n\}^n$ are {\sl formally conjugated} if there
exists an invertible $\Phi\in\C_0[[z_1,\ldots,z_n]]^n$ such that
$f_1=\Phi^{-1}\circ f_2\circ\Phi$ in~$\C_0[[z_1,\ldots,z_n]]^n$. 

It is clear that two holomorphically locally conjugated holomorphic local
dynamical systems are both formally and topologically locally conjugated too. On
the other hand, we shall see examples of holomorphic local dynamical systems
that are topologically locally conjugated without being neither formally nor
holomorphically locally conjugated, and examples of holomorphic local
dynamical systems that are formally conjugated without being
neither holomorphically nor topologically locally conjugated. So the last
natural question in local holomorphic dynamics we shall deal with is
\smallskip
\item{(Q5)}{\it Find normal forms and invariants with respect to the relation
of formal conjugacy for holomorphic local dynamical systems at~$O\in\C^n$.}
\smallskip
\noindent In this survey we shall present some of the main results known on
these questions, starting from the one-dimensional situation. 
But before entering the main core of the paper I would like to heartily thank Fran\c cois Berteloot, Kingshook Biswas, Filippo Bracci, Santiago Diaz-Madrigal, 
Graziano Gentili, Giorgio Patrizio, Mohamad Pouryayevali, Jasmin Raissy
and Francesca Tovena, without whom none of this would have been written.

\smallsect 2. One complex variable: the hyperbolic case

Let us then start by discussing holomorphic local dynamical systems at~$0\in\C$.
As remarked in the previous section, such a system
is given by a converging power series~$f$ without constant term:
$$
f(z)=a_1z+a_2z^2+a_3z^3+\cdots\in\C_0\{z\}\;.
$$

\Def The number $a_1=f'(0)$ is the {\sl multiplier} of~$f$.

Since $a_1 z$ is the best linear approximation of~$f$, it is sensible to expect
that the local dynamics of~$f$ will be strongly influenced by the value of~$a_1$.
For this reason we introduce the following definitions:

\Def Let $a_1\in\C$ be the multiplier of $f\in\End(\C,0)$. Then
\smallskip
\item{--} if $|a_1|<1$ we say that the fixed point $0$ is {\sl attracting;}
\item{--} if $a_1=0$ we say that the fixed point $0$ is {\sl superattracting;}
\item{--} if $|a_1|>1$ we say that the fixed point $0$ is {\sl repelling;}
\item{--} if $|a_1|\ne 0$,~1 we say that the fixed point $0$ is {\sl
hyperbolic;}
\item{--} if $a_1\in S^1$ is a root of unity, we say that the fixed point~$0$ is
{\sl parabolic} (or {\sl rationally indifferent\/});
\item{--} if $a_1\in S^1$ is not a root of unity, we say that the fixed 
point~$0$ is {\sl elliptic} (or {\sl irrationally indifferent\/}).

\noindent As we shall see in a minute, the dynamics of one-dimensional
holomorphic local dynamical systems with a hyperbolic fixed point is pretty
elementary; so we start with this case. 

\Rem Notice that if 0 is an attracting fixed
point for~$f\in\End(\C,0)$ with non-zero multiplier, then it is a repelling
fixed point for the inverse map~$f^{-1}\in\End(\C,0)$. 

Assume first that $0$ is attracting for the holomorphic local
dynamical system~$f\in\End(\C,0)$. Then we can write
$f(z)=a_1z+O(z^2)$, with~$0<|a_1|<1$; hence we can find a large constant~$M>0$,
a small constant~$\eps>0$ and
$0<\delta<1$ such that if $|z|<\eps$ then
$$
|f(z)|\le (|a_1|+M\eps)|z|\le\delta|z|\;.
\neweq\eqduuno
$$
In particular, if~$\Delta_\eps$ denotes the disk of center~$0$ and
radius~$\eps$, we have $f(\Delta_\eps)\subset\Delta_\eps$ for $\eps>0$ small
enough, and the stable set of~$f|_{\Delta_\eps}$ 
is~$\Delta_\eps$ itself (in particular, a one-dimensional attracting fixed point
is always stable). Furthermore,
$$
|f^k(z)|\le\delta^k|z|\to 0
$$
as $k\to+\infty$, and thus every orbit starting in~$\Delta_\eps$ is attracted by
the origin, which is the reason of the name ``attracting" for such a fixed
point. 

If instead 0 is a repelling fixed point, a similar argument (or the observation
that 0 is attracting for~$f^{-1}$) shows that for~$\eps>0$ small enough the
stable set of~$f|_{\Delta_\eps}$ reduces to the origin only: all (non-trivial)
orbits escape.

It is also not difficult to find holomorphic and topological normal forms for
one-dimensional holomorphic local dynamical systems with a hyperbolic fixed
point, as shown in the
following result, which can be considered as the beginning of the theory of holomorphic
dynamical systems:

\newthm Theorem \Koenigs: (K\oe nigs, 1884 [K\oe]) Let $f\in\End(\C,0)$ be a
one-dimensional holomorphic local dynamical system with a hyperbolic fixed
point at the origin, and let $a_1\in\C^*\setminus S^1$ be its multiplier. Then:
{\smallskip
\item{\rm(i)} $f$ is
holomorphically (and hence formally) locally conjugated to its linear
part~$g(z)=a_1 z$. The conjugation $\phe$ is uniquely determined by the
condition~$\phe'(0)=1$.
\item{\rm (ii)} Two such holomorphic local dynamical systems are
holomorphically conjugated if and only if they have the same multiplier.
\item{\rm (iii)} $f$ is topologically locally conjugated to the
map $g_<(z)=z/2$ if $|a_1|<1$, and to the map $g_>(z)=2z$\break
\indent if~$|a_1|>1$.}

\pf Let us assume $0<|a_1|<1$; if
$|a_1|>1$ it will suffice to apply the same argument to~$f^{-1}$.
\smallskip
(i) Choose $0<\delta<1$ such that $\delta^2<|a_1|<\delta$. Writing $f(z)=a_1
z+z^2r(z)$ for a suitable holomorphic germ~$r$, we can clearly find
$\eps>0$ such that $|a_1|+M\eps<\delta$, where
$M=\max_{z\in\bar{\Delta_\eps}}|r(z)|$. So we have
$$
|f(z)-a_1z|\le M|z|^2
$$
and
$$
|f^k(z)|\le\delta^k |z|
$$
for all $z\in\bar{\Delta_\eps}$ and $k\in\N$.

Put $\phe_k=f^k/a_1^k$; we claim that the
sequence $\{\phe_k\}$ converges to a holomorphic map
$\phe\colon\Delta_\eps\to\C$. Indeed we have
$$
|\phe_{k+1}(z)-\phe_k(z)|={1\over|a_1|^{k+1}}\bigl|f\bigl(f^k(z)\bigr)-a_1f^k(z)
\bigr|\le {M\over|a_1|^{k+1}}|f^k(z)|^2\le{M\over|a_1|}\left(
{\delta^2\over|a_1|}\right)^k|z|^2
$$
for all $z\in\bar{\Delta_\eps}$, and so the telescopic series
$\sum_k(\phe_{k+1}-\phe_k)$ is uniformly convergent in~$\Delta_\eps$
to~$\phe-\phe_0$.

Since $\phe_k'(0)=1$ for
all~$k\in\N$, we have $\phe'(0)=1$ and so, up to possibly shrink~$\eps$, we can
assume that $\phe$ is a biholomorphism with its image. Moreover, we have
$$
\phe\bigl(f(z)\bigr)=\lim_{k\to+\infty}{f^k\bigl(f(z)\bigr)\over a_1^k}
=a_1\lim_{k\to+\infty}{f^{k+1}(z)\over a_1^{k+1}}=a_1\phe(z)\;,
$$
that is $f=\phe^{-1}\circ g\circ \phe$, as claimed.

If $\psi$ is another local holomorphic function such that $\psi'(0)=1$
and $\psi^{-1}\circ g\circ\psi=f$, it follows that $\psi\circ\phe^{-1}(\lambda z)=
\lambda \psi\circ\phe^{-1}(z)$; comparing the expansion in power series of both
sides we find $\psi\circ\phe^{-1}\equiv\id$, that is $\psi\equiv\phe$, as claimed.

\smallskip
(ii) Since $f_1=\phe^{-1}\circ f_2\circ\phe$ implies $f_1'(0)=f_2'(0)$, the
multiplier is invariant under holomorphic local conjugation, and so two
one-dimensional holomorphic local dynamical systems with a hyperbolic fixed point
are holomorphically locally conjugated if and only if they have the same
multiplier.

\smallskip
(iii) Since $|a_1|<1$ it is easy to build a topological conjugacy between
$g$ and~$g_<$ on~$\Delta_\eps$. First choose a homeomorphism~$\chi$
between the annulus $\{|a_1|\eps\le
|z|\le\eps\}$ and the annulus $\{\eps/2\le |z|\le\eps\}$ which is the identity
on the outer circle and given by $\chi(z)=z/(2a_1)$ on the inner circle. Now
extend $\chi$ by induction to a homeomorphism between
the annuli $\{|a_1|^k\eps\le |z|\le|a_1|^{k-1}\eps\}$ and $\{\eps/2^k\le
|z|\le\eps/2^{k-1}\}$ by prescribing 
$$
\chi(a_1z)=\textstyle{{1\over2}}\chi(z)\;.
$$
Putting finally $\chi(0)=0$ we then get a homeomorphism $\chi$ of $\Delta_\eps$
with itself such that $g=\chi^{-1}\circ g_<\circ\chi$, as required.
\qedn

\Rem Notice that $g_<(z)={1\over 2}z$ and $g_>(z)=2z$ cannot be topologically
conjugated, because (for instance) $K_{g_<}$ is open whereas
$K_{g_>}=\{0\}$ is not.

\Rem The proof of this theorem is based on two techniques often used in dynamics to build conjugations. The first one is used in part (i). Suppose
that we would like to prove that two invertible local dynamical systems 
$f$, $g\in\End(M,p)$ are conjugated. Set $\phe_k=g^{-k}\circ f^k$, so that
$$
\phe_k\circ f=g^{-k}\circ f^{k+1}=g\circ\phe_{k+1}\;.
$$
Therefore if we can prove that $\{\phe_k\}$ converges to an invertible map 
$\phe$ as $k\to+\infty$ we get $\phe\circ f=g\circ\phe$, and thus $f$ and 
$g$ are conjugated, as desired. This is exactly the way we proved Theorem~\Koenigs.(i);
and we shall see variations of this techniques later on.

To describe the second technique we need a definition.

\Def Let $f\colon X\to X$ be an open continuous self-map of a topological 
space~$X$. A {\sl fundamental domain} for~$f$ is an open subset 
$D\subset X$ such that
\smallskip
\item{(i)} $f^h(D)\cap f^k(D)=\void$ for every $h\ne k\in\N$;
\item{(ii)} $\bigcup\limits_{k\in\N}f^k(\bar{D})=X$;
\item{(iii)} if $z_1$, $z_2\in\bar{D}$ are so that $f^h(z_1)=f^k(z_2)$
for some $h>k\in\N$ then $h=k+1$ and $z_2=f(z_1)\in\de D$.
\smallskip
\noindent There are other possible definitions of a fundamental domain,
but this will work for our aims. 

Suppose that we would like to prove that two open continuous maps $f_1\colon X_1
\to X_1$ and $f_2\colon X_2\to X_2$ are topologically conjugated. Assume we have fundamental domains $D_j\subset X_j$ for~$f_j$ (with $j=1$,~2) and a 
homeomorphism $\chi\colon\bar{D_1}\to\bar{D_2}$ such that
$$
\chi\circ f_1=f_2\circ\chi
\neweq\eqcfd
$$
on $\bar{D_1}\cap f_1^{-1}(\bar{D_1})$. Then we can extend $\chi$ to
a homeomorphism $\tilde\chi\colon X_1\to X_2$ conjugating $f_1$ and $f_2$
by setting
$$
\forevery{z\in X_1}
\tilde\chi(z)=f_2^{k}\bigl(\chi(w)\bigr)\;,
\neweq\equtd
$$ 
where $k=k(z)\in\N$ and $w=w(z)\in\bar{D}$ are chosen so that $f_1^k(w)=z$. 
The definition of fundamental domain and \eqcfd\ imply that $\tilde\chi$
is well-defined. Clearly $\tilde\chi\circ f_1=f_2\circ\tilde\chi$; and 
using the openness of $f_1$ and $f_2$ it is easy to
check that $\tilde\chi$ is a homeomorphism. This is the technique we used
in the proof of Theorem~\Koenigs.(iii); and we shall use it again later.

Thus the dynamics in the one-dimensional hyperbolic case is completely clear.
The superattracting case can be treated similarly. If $0$ is a superattracting
point for an~$f\in\End(\C,0)$, we can write 
$$
f(z)=a_rz^r+a_{r+1}z^{r+1}+\cdots
$$
with $a_r\ne 0$. 

\Def The number $r\ge 2$ is the {\sl order} (or {\sl local degree\/}) of the superattracting
point. 

An argument similar to the one described before shows that for $\eps>0$
small enough the stable set of~$f|_{\Delta_\eps}$ still is all of~$\Delta_\eps$,
and the orbits converge (faster than in the attracting case) to the origin.
Furthermore, we can prove the following

\newthm Theorem \Bottcher: (B\"ottcher, 1904 [B\"o]) Let $f\in\End(\C,0)$ be a
one-dimensional holomorphic local dynamical system with a superattracting fixed
point at the origin, and let $r\ge 2$ be its order. Then:
{\smallskip
\item{\rm(i)} $f$ is
holomorphically (and hence formally) locally conjugated to the map~$g(z)=z^r$, and the conjugation is unique up to multiplication by an $(r-1)$-root
of unity;
\itm{(ii)} two such holomorphic local dynamical systems are holomorphically
(or topologically) conjugated if and\break\indent only if they have the same
order.}

\pf First of all, up to a linear conjugation $z\mapsto\mu z$ with
$\mu^{r-1}=a_r$ we can assume $a_r=1$.

Now write $f(z)=z^r h_1(z)$ for a suitable holomorphic germ $h_1$ with
$h_1(0)=1$. By induction, it is easy to see that we can write
$f^k(z)=z^{r^k}h_k(z)$ for a suitable holomorphic germ $h_k$ with $h_k(0)=1$.
Furthermore, the equalities $f\circ f^{k-1}=f^k=f^{k-1}\circ f$ yields
$$
h_{k-1}(z)^r
h_1\bigl(f^{k-1}(z)\bigr)=h_k(z)=h_1(z)^{r^{k-1}}h_{k-1}\bigl(f(z)
\bigr)\;.
\neweq\eqddue
$$
Choose $0<\delta<1$. Then we can clearly find $1>\eps>0$ such that
$M\eps<\delta$, where $M=\max_{z\in\bar{\Delta_\eps}}|h_1(z)|$; we can also
assume that
$h_1(z)\ne 0$ for all $z\in\bar{\Delta_\eps}$. Since
$$
\forevery{z\in\bar{\Delta_\eps}}|f(z)|\le M|z|^r<\delta |z|^{r-1}\;,
$$
we have $f(\Delta_\eps)\subset\Delta_\eps$, as anticipated before. 

We also remark that \eqddue\ implies that each $h_k$ is well-defined and never
vanishing on~$\bar{\Delta_\eps}$.
So for every~$k\ge 1$ we can choose a unique $\psi_k$ holomorphic in
$\Delta_\eps$ such that $\psi_k(z)^{r^k}=h_k(z)$ on $\Delta_\eps$ and
with~$\psi_k(0)=1$. 

Set $\phe_k(z)=z\psi_k(z)$, so that $\phe'_k(0)=1$ and
$\phe_k(z)^{r^k}=f_k(z)$ on $\Delta_\eps$; in particular, formally we have
$\phe_k=g^{-k}\circ f^k$. We claim that the sequence
$\{\phe_k\}$ converges to a holomorphic function $\phe$ on $\Delta_\eps$.
Indeed, we have
$$
\eqalign{
\left|{\phe_{k+1}(z)\over\phe_k(z)}\right|&=\left|{\psi_{k+1}(z)^{r^{k+1}}\over
\psi_k(z)^{r^{k+1}}}\right|^{1/r^{k+1}}=\left|{h_{k+1}(z)\over h_k(z)^r}\right|
^{1/r^{k+1}}=\bigl|h_1\bigl(f^k(z)\bigr)\bigr|^{1/r^{k+1}}\cr
&=\bigl|1+O\bigl(|f^k(z)|\bigr)\bigr|^{1/r^{k+1}}=1+{1\over
r^{k+1}}O\bigl(|f^k(z)|\bigr)=1+O\left({1\over r^{k+1}}\right)\;,
\cr}
$$
and so the telescopic product $\prod_k(\phe_{k+1}/\phe_k)$ converges to
$\phe/\phe_1$ uniformly in $\Delta_\eps$.

Since $\phe_k'(0)=1$ for
all~$k\in\N$, we have $\phe'(0)=1$ and so, up to possibly shrink~$\eps$, we can
assume that $\phe$ is a biholomorphism with its image. Moreover, we have
$$
\phe_k\bigl(f(z)\bigr)^{r^k}=f(z)^{r^k}\psi_k\bigl(f(z)\bigr)^{r^k}=
z^{r^k}h_1(z)^{r^k}h_k\bigl(f(z)\bigr)=z^{r^{k+1}}h_{k+1}(z)=
\bigl[\phe_{k+1}(z)^r\bigr]^{r^k}\;,
$$
and thus $\phe_k\circ f=[\phe_{k+1}]^r$. Passing to the limit we get
$f=\phe^{-1}\circ g\circ \phe$, as claimed. 

If $\psi$ is another local biholomorphism conjugating $f$ with $g$,
we must have $\psi\circ\phe^{-1}(z^r)=\psi\circ\phe^{-1}(z)^r$ for
all $z$ in a neighbourhood of the origin; comparing the series 
expansions at the origin we get $\psi\circ\phe^{-1}(z)=az$ with 
$a^{r-1}=1$, and hence $\psi(z)=a\phe(z)$, as claimed.

Finally, (ii) follows because $z^r$ and $z^s$ are locally
topologically conjugated if and only if~$r=s$ (because the order is
the number of preimages of points close to the origin).\qedn

Therefore the one-dimensional local dynamics about a hyperbolic or
superattracting fixed point is completely clear; let us now discuss what happens
about a parabolic fixed point.

\smallsect 3. One complex variable: the parabolic case

Let $f\in\End(\C,0)$ be a (non-linear) holomorphic local dynamical system with a
parabolic fixed point at the origin. Then we can write 
$$
f(z)=e^{2i\pi p/q}z+a_{r+1} z^{r+1}+a_{r+2}z^{r+2}+\cdots,
\neweq\eqtuno
$$
with $a_{r+1}\ne 0$.

\Def The rational number $p/q\in\Q\cap[0,1)$ is the {\sl rotation number}
of~$f$, and the number $r+1\ge 2$ is the {\sl multiplicity} of~$f$ at the 
fixed point. If $p/q=0$ (that is, if the multiplier is~1), we shall say
that $f$ is {\sl tangent to the identity.}

The first observation is that such a dynamical system is never locally
conjugated to its linear part, not even topologically, unless it is of 
finite order:

\newthm Proposition \tunoa: Let $f\in\End(\C,0)$ be a holomorphic local
dynamical system with multiplier $\lambda$, and assume that $\lambda=e^{2i\pi p/q}$ is a
primitive root of the unity of order $q$. Then $f$ is holomorphically (or
topologically or formally) locally conjugated to $g(z)=\lambda z$ if and only if $f^q\equiv\id$.

\pf If $\phe^{-1}\circ f\circ\phe(z)=e^{2\pi ip/q}z$ then $\phe^{-1}\circ
f^q\circ\phe=\id$, and hence~$f^q=\id$.

Conversely, assume that $f^q\equiv\id$ and set
$$
\phe(z)={1\over q}\sum_{j=0}^{q-1}{f^j(z)\over\lambda^j}.
$$
Then it is easy to check that $\phe'(0)=1$ and $\phe\circ f(z)=\lambda\phe(z)$,
and so $f$ is holomorphically (and topologically and formally) 
locally conjugated to~$\lambda z$.\qedn

In particular, if $f$ is tangent to the identity then it {\it
cannot} be locally conjugated to the identity (unless it was the identity
to begin with, which is not a very interesting case dynamically
speaking). More precisely, the stable set of such an
$f$ is never a neighbourhood of the origin. To understand why, let us first
consider a map of the form
$$
f(z)=z(1+az^r)
$$
for some $a\ne 0$. Let $v\in S^1\subset\C$ be such that $av^r$ is real and
positive. Then for any $c>0$ we have
$$
f(cv)=c(1+c^rav^r)v\in\R^+v;
$$
moreover, $|f(cv)|>|cv|$. In other words, the half-line~$\R^+v$ is $f$-invariant
and repelled from the origin, that is $K_f\cap\R^+ v=\void$. Conversely, if
$av^r$ is real and negative then the segment $[0,|a|^{-1/r}]v$ is $f$-invariant
and attracted by the origin. So $K_f$ neither is a neighbourhood of the origin
nor reduces to~$\{0\}$. 

This example suggests the following definition:

\Def Let $f\in\End(\C,0)$ be tangent to the identity of multiplicity~$r+1\ge2$. Then a unit vector~$v\in S^1$
is an {\sl attracting} (respectively, {\sl repelling\/}) {\sl direction}
for~$f$ at the origin if~$a_{r+1}v^r$ is real and negative (respectively,
positive). 

Clearly, there are $r$ equally spaced attracting directions,
separated by $r$ equally spaced repelling directions: if
$a_{r+1}=|a_{r+1}|e^{i\alpha}$, then $v=e^{i\theta}$ is attracting
(respectively, repelling) if and only if
$$
\theta={2k+1\over r}\pi-{\alpha\over
r}\qquad\Biggl(\hbox{respectively, }\theta={2k\over r}\pi-{\alpha\over r}
\Biggr).
$$
Furthermore, a repelling
(attracting) direction for~$f$ is attracting (repelling) for~$f^{-1}$, which is
defined in a neighbourhood of the origin.

It turns out that to every attracting direction is associated a connected
component of~$K_f\setminus\{0\}$. 

\Def Let~$v\in S^1$ be an attracting direction for
an~$f\in\End(\C,0)$ tangent to the identity. The {\sl basin} centered at~$v$ is the set of
points~$z\in K_f\setminus\{0\}$ such that $f^k(z)\to 0$ and $f^k(z)/|f^k(z)|\to
v$ (notice that, up to shrinking the domain of~$f$, we can assume that $f(z)\ne
0$ for all $z\in K_f\setminus\{0\}$). If $z$ belongs to the basin centered
at~$v$, we shall say that the orbit of~$z$ {\sl tends to~$0$ tangent to~$v$.}

A slightly more specialized (but more useful) object is the following: 

\Def An {\sl
attracting petal} centered at an attracting direction~$v$ of
an~$f\in\End(\C,0)$ tangent to the identity is an open simply
connected $f$-invariant set $P\subseteq K_f\setminus\{0\}$ such that a point
$z\in K_f\setminus\{0\}$ belongs to the basin centered at~$v$ if and only if its
orbit intersects~$P$. In other words, the orbit of a point tends to~$0$ tangent
to~$v$ if and only if it is eventually contained in~$P$. A {\sl repelling petal}
(centered at a repelling direction) is an attracting petal for the inverse
of~$f$.

It turns out that the basins centered at the attracting directions are exactly
the connected components of~$K_f\setminus\{0\}$, as shown in the {\it Leau-Fatou
flower theorem:} 

\newthm Theorem \flower: (Leau, 1897 [L]; Fatou, 1919-20 [F1--3]) Let
$f\in\End(\C,0)$ be a holomorphic local dynamical system tangent to the identity
with multiplicity~$r+1\ge 2$ at the fixed point. Let $v_1^+,\ldots,v_r^+\in
S^1$ be the
$r$ attracting directions of~$f$ at the origin, and $v_1^-,\ldots,v_r^-\in
S^1$ the $r$ repelling directions. Then {\smallskip
\item{\rm (i)} for each attracting (repelling) direction~$v_j^\pm$
there exists an attracting (repelling) petal~$P_j^\pm$, so that the
union of these $2r$ petals together with the origin forms a neighbourhood
of the origin. Furthermore, the
$2r$ petals are arranged ciclically so that two petals intersect if and only if
the angle between their central directions is~$\pi/r$.
\item{\rm (ii)} $K_f\setminus\{0\}$ is the (disjoint) union of the basins
centered at the
$r$ attracting directions.
\item{\rm (iii)} If $B$ is a basin centered at one of the attracting directions, then
there is a function $\phe\colon B\to\C$ such that\break\indent $\phe\circ
f(z)=\phe(z)+1$ for all $z\in B$. Furthermore, if $P$ is the corresponding petal
constructed in part {\rm (i),}\break\indent then $\phe|_P$ is a biholomorphism
with an open subset of the complex plane containing a right
half-plane\break\indent --- and so $f|_P$ is holomorphically
conjugated to the translation $z\mapsto z+1$.}

\pf Up to a linear conjugation, we
can assume that $a_{r+1}=-1$, so that the attracting directions are the $r$-th
roots of unity. For any $\delta>0$, the set $\{z\in\C\mid|z^r-\delta|<\delta\}$
has exactly $r$ connected components, each one symmetric with respect to a
different $r$-th root of unity; it will turn out that, for $\delta$ small
enough, these connected components are attracting petals of~$f$, even though to
get a pointed neighbourhood of the origin we shall need larger petals.

For $j=1,\ldots,r$ let $\Sigma_j\subset\C^*$ denote the sector centered
about the attractive direction~$v_j^+$ and bounded by two consecutive repelling
directions, that is
$$
\Sigma_j=\left\{z\in\C^*\biggm| {2j-3\over r}\,\pi<\arg z<{2j-1\over
r}\,\pi\right\}\;.
$$
Notice that each $\Sigma_j$ contains a unique connected component $P_{j,\delta}$
of $\{z\in\C\mid|z^r-\delta|<\delta\}$; moreover, $P_{j,\delta}$ is tangent
at the origin to the sector centered about~$v_j$ of amplitude~$\pi/r$. 

The main technical trick in this proof consists in transfering the setting to a
neighbourhood of infinity in the Riemann sphere~$\P^1(\C)$. Let
$\psi\colon\C^*\to\C^*$ be given by 
$$
\psi(z)={1\over rz^r}\;;
$$
it is a biholomorphism between $\Sigma_j$ and
$\C^*\setminus\R^-$, with inverse $\psi^{-1}(w)=(rw)^{-1/r}$, choosing
suitably the $r$-th root. Furthermore, $\psi(P_{j,\delta})$ is
the right half-plane
$H_\delta=\{w\in\C\mid\Re w>1/(2r\delta)\}$.

When $|w|$ is so large that $\psi^{-1}(w)$ belongs to the domain of definition
of~$f$, the composition $F=\psi\circ f\circ\psi^{-1}$ makes sense, and we have
$$
F(w)=w+1+O(w^{-1/r})\;.
\neweq\eqtdue
$$
Thus to study the dynamics of $f$ in a neighbourhood of the origin in~$\Sigma_j$
it suffices to study the dynamics of~$F$ in a neighbourhood of infinity.

The first observation is that when $\Re w$ is large enough then
$$
\Re F(w)>\Re w+{1\over 2}\;;
$$  
this implies that for $\delta$ small enough $H_\delta$ is $F$-invariant (and
thus $P_{j,\delta}$ is $f$-invariant). Furthermore, by induction one has
$$
\forevery{w\in H_\delta} \Re F^k(w)>\Re w+{k\over 2}\;,
\neweq\eqttre
$$
which implies that $F^k(w)\to\infty$ in $H_\delta$ (and $f^k(z)\to 0$ in
$P_{j,\delta}$) as $k\to\infty$.

Now we claim that the argument of $w_k=F^k(w)$ tends to zero. Indeed, \eqtdue\
and \eqttre\ yield
$$
{w_k\over k}={w\over k}+1+{1\over k}\sum_{l=0}^{k-1} O(w_l^{-1/r})\;;
$$
so Cesaro's theorem on the averages of a converging sequence implies 
$$
{w_k\over k}\to1\;,
\neweq\eqtquattroa
$$
and thus $\arg w_k\to 0$ as~$k\to\infty$.  
Going back to~$P_{j,\delta}$, this implies that $f^k(z)/|f^k(z)|\to v_j$ for
every $z\in P_{j,\delta}$. Since furthermore $P_{j,\delta}$ is centered
about~$v_j^+$, every orbit converging to~$0$ tangent to~$v_j^+$ must
intersect~$P_{j,\delta}$, and thus we have proved that $P_{j,\delta}$ is an
attracting petal.

Arguing in the same way with~$f^{-1}$ we get repelling petals; unfortunately, the petals obtained so far
are too small to form a full pointed neighbourhood of the origin. In fact, as
remarked before each
$P_{j,\delta}$ is contained in a sector centered about~$v_j$ of amplitude~$\pi/r$; therefore
the repelling and attracting petals obtained in this way do not intersect but are tangent to
each other.  We need larger petals.

So our aim is to find an $f$-invariant subset $P_j^+$ of $\Sigma_j$
containing~$P_{j,\delta}$ and which is tangent at the origin to a sector centered
about~$v_j^+$ of amplitude strictly greater than~$\pi/r$. To do so, first of all
remark that there are $R$, $C>0$ such that
$$
|F(w)-w-1|\le {C\over|w|^{1/r}}
\neweq\eqtcinquea
$$
as soon as $|w|>R$. Choose
$\eps\in(0,1)$ and select~$\delta>0$ so that $4r\delta<R^{-1}$ and 
$\eps>2C(4r\delta)^{1/r}$. Then
$|w|>1/(4r\delta)$ implies
$$
|F(w)-w-1|< \eps/2\;.
$$
Set $M_\eps=(1+\eps)/(2r\delta)$ and let
$$
\tilde H_\eps=\{w\in\C\mid |\Im w|>-\eps\Re w+M_\eps\}\cup H_\delta\;.
$$
If $w\in\tilde H_\eps$ we have $|w|>1/(2r\delta)$ and hence
$$
\Re F(w)>\Re w +1-\eps/2\qquad\hbox{and}\qquad |\Im F(w)-\Im w|<\eps/2\;;
\neweq\eqbht
$$
it is then easy to check that $F(\tilde H_\eps)\subset\tilde H_\eps$ and
that every orbit starting in~$\tilde H_\eps$ must eventually enter~$H_\delta$.
Thus $P_j^+=\psi^{-1}(\tilde H_\eps)$ is as required, and we have proved
(i).

To prove (ii) we need a further property of $\tilde H_\eps$. If $w\in\tilde H_\eps$, arguing by induction on $k\ge 1$ using \eqbht\ we get
$$
k\left(1-{\eps\over2}\right)<\Re F^k(w)-\Re w
$$
and
$$
{k\eps(1-\eps)\over 2}<
|\Im F^k(w)|+\eps\Re F^k(w)-\bigl(|\Im w|+\eps\Re w\bigr)\;.
$$
This implies that for every $w_0\in\tilde H_\eps$ there exists a $k_0\ge 1$
so that we cannot have $F^{k_0}(w)=w_0$ for any~$w\in\tilde H_\eps$.
Coming back to the $z$-plane, this says that any inverse orbit of $f$
must eventually leave~$P_j^+$. Thus every (forward) orbit of
$f$ must eventually leave any repelling petal. So if $z\in K_f\setminus\{O\}$,
where the stable set is computed working in the neighborhood of the origin
given by the union of repelling and attracting petals (together with
the origin), the orbit of $z$ must eventually land in an attracting
petal, and thus $z$ belongs to a basin centered at one of the $r$ attracting
directions --- and (ii) is proved.
  
To prove (iii), first of all we notice that we have
$$
|F'(w)-1|\le {2^{1+1/r}C\over|w|^{1+1/r}}
\neweq\eqtseia
$$
in $\tilde H_\eps$. Indeed, \eqtcinquea\ says that if $|w|>1/(2r\delta)$ then
the function $w\mapsto F(w)-w-1$ sends the disk of center $w$ and radius $|w|/2$
into the disk of center the origin and radius $C/(|w|/2)^{1/r}$; inequality
\eqtseia\ then follows from the Cauchy estimates on the derivative.

Now choose $w_0\in H_\delta$, and set
$\tilde\phe_k(w)=F^k(w)-F^k(w_0)$. Given $w\in\tilde H_\eps$, as soon as
$k\in\N$ is so large that $F^k(w)\in H_\delta$ we can apply Lagrange's theorem to
the segment from $F^k(w_0)$ to $F^k(w)$ to get a $t_k\in[0,1]$ such that
$$
\eqalign{
\left|{\tilde\phe_{k+1}(w)\over\tilde\phe_k(w)}-1\right|&=
\left|{F\bigl(F^k(w)\bigr)-F^k\bigl(F^k(w_0)\bigr)\over F^k(w)-F^k(w_0)}-1\right|
=\bigl|F'\bigl(t_kF^k(w)+(1-t_k)F^k(w_0)\bigr)-1\bigr|\cr
&\le {2^{1+1/r}C\over\min\{\Re|F^k(w)|,\Re|F^k(w_0)|\}^{1+1/r}}\le
	{C'\over k^{1+1/r}}\;,
\cr}
$$
where we used \eqtseia\ and \eqtquattroa, and the constant $C'$ is uniform on
compact subsets of $\tilde H_\eps$ (and it can be chosen uniform on $H_\delta$).

As a consequence, the telescopic product $\prod_k\tilde\phe_{k+1}/\tilde\phe_k$
converges uniformly on compact subsets of $\tilde H_\eps$ (and uniformly on
$H_\delta$), and thus the sequence
$\tilde\phe_k$ converges, uniformly on compact subsets, to a holomorphic
function $\tilde\phe\colon\tilde H_\eps\to\C$. Since we have
$$
\tilde\phe_k\circ F(w)=F^{k+1}(w)-F^k(w_0)=\tilde\phe_{k+1}(w)+F\bigl(F^k(w_0)
\bigr)-F^k(w_0)=\tilde\phe_{k+1}(w)+1+O\bigl(|F^k(w_0)|^{-1/r}\bigr)\;,
$$
it follows that
$$
\tilde\phe\circ F(w)=\tilde\phe(w)+1
$$
on $\tilde H_\eps$. In particular, $\tilde\phe$ is not constant;
being the limit of injective functions, by Hurwitz's theorem it is injective.

We now prove that the image of $\tilde\phe$ contains a right half-plane. First
of all, we claim that
$$
\lim_{{|w|\to+\infty}\atop w\in H_\delta}{\tilde\phe(w)\over w}=1\;.
\neweq\eqtsette
$$
Indeed, choose $\eta>0$. Since the convergence of the telescopic product is
uniform on~$H_\delta$, we can find $k_0\in\N$ such that
$$
\left|{\tilde\phe(w)-\tilde\phe_{k_0}(w)\over w-w_0}\right|<{\eta\over3}
$$
on $H_\delta$. 
Furthermore, we have
$$
\left|{\tilde\phe_{k_0}(w)\over w-w_0}-1\right|=\left|{k_0+\sum_{j=0}^{k_0-1}
O(|F^j(w)|^{-1/r})+w_0-F^{k_0}(w_0)\over w-w_0}\right|=O(|w|^{-1})
$$
on $H_\delta$; therefore we can find $R>0$ such that 
$$
\left|{\tilde\phe(w)\over w-w_0}-1\right|<{\eta\over3}
$$
as soon as $|w|>R$ in $H_\delta$. Finally, if $R$ is large enough we also 
have
$$
\left|{\tilde\phe(w)\over w-w_0}-{\tilde\phe(w)\over w}\right|=
\left|{\tilde\phe(w)\over w-w_0}\right|\left|{w\over w_0}\right|
<{\eta\over3}\;,
$$
and \eqtsette\ follows.

Equality \eqtsette\ clearly implies that $(\tilde\phe(w)-w^o)/(w-w^o)\to 1$ as
$|w|\to+\infty$ in~$H_\delta$ for any $w^o\in\C$. But this means that if $\Re
w^o$ is large enough then the difference between the variation of the argument of
$\tilde\phe-w^o$ along a suitably small closed circle around $w^o$ and
the variation of the argument of $w-w^o$ along the same circle will be less than
$2\pi$ --- and thus it will be zero. Then the argument principle implies
that $\tilde\phe-w^o$ and $w-w^o$ have the same number of zeroes inside that
circle, and thus $w^o\in\tilde\phe(H_\delta)$, as required.

So setting $\phe=\tilde\phe\circ\psi$, we have defined a function $\phe$ with
the required properties on $P_j^+$. To extend it to the whole basin $B$ it
suffices to put
$$
\phe(z)=\phe\bigl(f^k(z)\bigr)-k\;,
\neweq\equtd
$$ 
where $k\in\N$ is the first integer such that $f^k(z)\in P_j^+$. \qedn

\Rem It is possible to construct petals that cannot be contained in any sector
strictly smaller than $\Sigma_j$. To do so we need an $F$-invariant
subset $\hat H_\eps$ of $C^*\setminus\R^-$ containing 
$\tilde H_\eps$  and containing eventually every half-line issuing from the
origin (but~$\R^-$). For $M>>1$ and $C>0$ large enough, replace the straight
lines bounding $\tilde H_\eps$ on the left of $\Re w=-M$ by the curves
$$
|\Im w|=\cases{C\log|\Re w|& if $r=1$,\cr
\noalign{\smallskip}
C|\Re w|^{1-1/r}& if $r>1$.
\cr}
$$
Then it is not too difficult to check that the domain $\hat H_\eps$ so obtained
is as desired (see [CG]). 


So we have a complete description of the dynamics in the neighbourhood of the
origin. Actually, Camacho has pushed this argument even further, obtaining a
complete topological classification of one-dimensional holomorphic local
dynamical systems tangent to the identity (see also [BH, Theorem~1.7]):

\newthm Theorem \Camacho: (Camacho, 1978 [C]; Shcherbakov, 1982 [S]) Let
$f\in\End(\C,0)$ be a holomorphic local dynamical system tangent to the identity
with multiplicity~$r+1$ at the fixed point. Then $f$ is topologically locally
conjugated to the map
$$
g(z)=z-z^{r+1}\;.
$$

The formal classification is simple too, though different 
(see, e.g., Milnor [Mi]):

\newthm Proposition \formaltangent: Let $f\in\End(\C,0)$ be a
holomorphic local dynamical system tangent to the identity with
multiplicity~$r+1$ at the fixed point. Then $f$ is formally conjugated to the map
$$
g(z)=z-z^{r+1}+\beta z^{2r+1}\;,
\neweq\eqtottoa
$$
where $\beta$ is a formal (and holomorphic) invariant given
by
$$
\beta={1\over 2\pi i}\int_\gamma{dz\over z-f(z)}\;,
\neweq\eqttre
$$
where the integral is taken over a small positive loop~$\gamma$ about the
origin.

\pf An easy computation shows that if $f$ is given by \eqtottoa\ then
\eqttre\ holds. Let us now show that the integral in \eqttre\ is a
holomorphic invariant. Let
$\phe$ be a local biholomorphism fixing the origin, and set $F=\phe^{-1}\circ
f\circ\phe$. Then
$$
{1\over2\pi i}\int_\gamma {dz\over z-f(z)}={1\over2\pi
i}\int_{\phe^{-1}\circ\gamma} {\phe'(w)\,dw\over \phe(w)-f\bigl(\phe(w)\bigr)}
={1\over2\pi
i}\int_{\phe^{-1}\circ\gamma} {\phe'(w)\,dw\over \phe(w)-\phe\bigl(F(w)\bigr)}\;.
$$
Now, we can clearly find $M$, $M_1>0$ such that
$$
\eqalign{
\left|{1\over w-F(w)}-{\phe'(w)\over\phe(w)-\phe\bigl(F(w)\bigr)}\right|&=
{1\over\bigl|\phe(w)-\phe\bigl(F(w)\bigr)\bigr|}
\left|{\phe(w)-\phe\bigl(F(w)\bigr)\over
w-F(w)}-\phe'(w)\right|\cr
&\le M{|w-F(w)|\over\bigl|\phe(w)-\phe\bigl(F(w)\bigr)\bigr|}\le M_1\;,
\cr}
$$
in a neighbourhood of the origin, where the last inequality follows from the fact
that $\phe'(0)\ne 0$. This means that the two meromorphic functions
$1/\bigl(w-F(w)\bigr)$ and $\phe'(w)/\bigl(\phe(w)-\phe(\bigl(F(w)\bigr)\bigr)$
differ by a holomorphic function; so they have the same integral along any small
loop surrounding the origin, and 
$$
{1\over2\pi i}\int_\gamma {dz\over z-f(z)}={1\over2\pi
i}\int_{\phe^{-1}\circ\gamma} {dw\over w-F(w)},
$$
as claimed.

To prove that $f$ is formally conjugated to $g$, let us first take a local
formal change of coordinates $\phe$ of the form
$$
\phe(z)=z+\mu z^d+O_{d+1}
\neweq\eqtuzero
$$
with $\mu\ne0$, and where we are writing $O_{d+1}$ instead of $O(z^{d+1})$. It
follows that $\phe^{-1}(z)=z-\mu z^d+O_{d+1}$, $(\phe^{-1})'(z)=1-d\mu
z^{d-1}+O_d$ and $(\phe^{-1})^{(j)}=O_{d-j}$ for all $j\ge 2$. Then using
the Taylor expansion of $\phe^{-1}$ we get
$$
\eqalign{
\phe^{-1}\circ f\circ\phe(z)&=\phe^{-1}\left(\phe(z)+\sum_{j\ge r+1}a_j\phe(z)^j
\right)\cr
&=z+(\phe^{-1})'\bigl(\phe(z)\bigr)\sum_{j\ge r+1} a_j z^j(1+\mu z^{d-1}+O_d)^j
+O_{d+2r}\cr
&=z+[1-d\mu z^{d-1}+O_d]\sum_{j\ge r+1}a_jz^j(1+j\mu z^{d-1}+O_d)+O_{d+2r}\cr
&=z+a_{r+1}z^{r+1}+\cdots+a_{r+d-1}z^{r+d-1}+[a_{r+d}+(r+1-d)\mu a_{r+1}] z^{r+d}
+O_{r+d+1}.
\cr}
\neweq\eqtuuno
$$
This means that if $d\ne r+1$ we can use a polynomial change of coordinates of
the form $\phe(z)=z+\mu z^d$ to remove the term of degree $r+d$ from the Taylor
expansion of $f$ without changing the lower degree terms. 

So to conjugate $f$ to $g$ it suffices to use a linear change of coordinates to
get $a_{r+1}=-1$, and then apply a sequence of change of coordinates of the form
$\phe(z)=z+\mu z^d$ to kill all the terms in the Taylor expansion of $f$ but
the term of degree $z^{2r+1}$. 

Finally, formula \eqtuuno\ also shows that two maps of the form \eqtottoa\
with different $\beta$ cannot be formally conjugated, and we are done.\qedn

\Def The number~$\beta$ given by~\eqttre\ is called {\sl index} of~$f$ at the fixed point. The {\sl iterative residue} of~$f$ is then defined by
$$
\hbox{Resit}(f)={r+1\over 2}-\beta\;.
$$

The iterative residue has been introduced by \'Ecalle [\'E1], and it behaves nicely under iteration; for instance, it is possible to prove (see [BH, Proposition~3.10]) that
$$
\hbox{Resit}(f^k)={1\over k}\,\hbox{Resit}(f)\;.
$$

The holomorphic classification of maps tangent to the identity is much more complicated: as shown by
\'Ecalle [\'E2--3] and Voronin [V] in~1981, it depends on functional invariants. 
We shall now try and
roughly describe it; see~[I2], [M1-2], [K], [BH] and the original
papers for details. 

Let $f\in\End(\C,0)$ be tangent to the
identity with multiplicity~$r+1$ at the fixed point; up to a linear
change of coordinates we can assume that~$a_{r+1}=-1$. Let
$P_j^\pm$ be a set of petals as in Theorem~\flower.(i), ordered so
that $P_1^+$ is centered on the positive real semiaxis, and the others are
arranged cyclically counterclockwise. Denote by
$\phe_j^+$ (respectively, $\phe_j^-$) the biholomorphism conjugating $f|_{P_j^+}$ (respectively, $f|_{P_j^-}$) to the shift~$z\mapsto z+1$ in a right (respectively, left) half-plane given by
Theorem~\flower.(iii) --- applied to~$f^{-1}$ for the repelling petals. If we
moreover require that 
$$
\phe_j^\pm(z)={1\over rz^r}\pm\hbox{Resit}(f)\cdot\log z+o(1)\;,
\neweq\eqtquattro
$$
then $\phe_j$ is uniquely determined. 

Put now $U_j^+=P_j^-\cap P_{j+1}^+$, $U_j^-=P_j^-\cap P_j^+$, and
$S_j^\pm=\bigcup_{k\in\Z}U_j^\pm$. Using the dynamics as in \equtd\ we
can extend $\phe_j^-$ to $S_j^\pm$, and $\phe_j^+$ to
$S^+_{j-1}\cup S^-_j$; put $V_j^\pm=\phe_j^-(S_j^\pm)$, $W_j^-=\phe_j^+
(S_j^-)$ and $W_j^+=\phe_{j+1}^+(S_j^+)$. Then let $H_j^-\colon V_j^-\to
W_J^-$ the restriction of~$\phe_j^+\circ(\phe_j^-)^{-1}$ to~$V_j^-$,
and $H_j^+\colon V_j^+\to W_j^+$ the restriction of $\phe_{j+1}^+\circ
(\phe_j^-)^{-1}$ to~$V_j^+$.

It is not difficult to see that $V_j^\pm$ and $W_j^\pm$ are invariant by translation by~1, and that $V_j^+$ and $W_j^+$ contain an upper
half-plane while $V_j^-$ and $W_j^-$ contain a lower half-plane.
Moreover, we have $H_j^\pm(z+1)=H_j^\pm(z)+1$; therefore using 
the projection $\pi(z)=\exp(2\pi iz)$ we can induce holomorphic
maps $h_j^\pm\colon\pi(V_j^\pm)\to\pi(W_j^\pm)$, 
where $\pi(V_j^+)$ and $\pi(W_j^+)$ are pointed neighbourhood of the
origin, and $\pi(V_j^-)$ and $\pi(W_j^-)$ are pointed neighbourhood
of $\infty\in\P^1(\C)$. 

It is possible to show that setting $h_j^+(0)=0$ one obtains a holomorphic germ $h_j^+\in\End(\C,0)$, and that setting $h_j^-(\infty)=\infty$ one obtains a holomorphic germ $h_j^+\in\End\bigl(\P^1(C),\infty\bigr)$.
Furthermore, denoting by $\lambda_j^+$ (respectively, $\lambda_j^-$) the multiplier of $h_j^+$ at~0 (respectively, of $h_j^-$ at~$\infty$), it turns out that
$$
\prod_{j=1}^r(\lambda_j^+\lambda_j^-)=\exp\bigr[4\pi^2\hbox{Resit}(f)
\bigr]\;.
\neweq\eqocon
$$
  
Now, if we replace $f$ by a holomorphic local conjugate~$\tilde f=\psi^{-1}\circ f\circ
\psi$, and denote by $\tilde h_j^\pm$ the corresponding germs, it is not difficult to check that (up to a cyclic renumbering of the petals)
there are constants $\alpha_j$, $\beta_j\in\C^*$ such that
$$
\tilde h_j^-(z)=\alpha_j h_j^-\left({z\over\beta_j}\right)\qquad\hbox{and}
\qquad \tilde h_j^+(z)=\alpha_{j+1} h^+_j\left({z\over\beta_j}\right)\;.
\neweq\eqdcon
$$
This
suggests the introduction of an equivalence relation on the set of $2r$-uple of
germs $(h_1^\pm,\ldots,h_r^\pm)$.

\Def Let $M_r$ denote the set of $2r$-uple of germs~${\bf h}=(h_1^\pm,\ldots,h_r^\pm)$, with $h_j^+\in\End(\C,0)$, $h_j^-\in\End\bigl(\P^1(\C),\infty\bigr)$, and whose multipliers satisfy \eqocon.
We shall say that ${\bf h}$,~$\tilde{\bf h}\in
M_r$ are {\sl equivalent} if up to a cyclic permutation of the indeces 
we have \eqdcon\ for suitable $\alpha_j$, $\beta_j\in\C^*$. We denote by~$\ca M_r$ the set of all equivalence classes.

The procedure described above allows then to associate to any $f\in\End(\C,0)$
tangent to the identity with multiplicity~$r+1$ an
element~$\mu_f\in\ca M_r$.

\Def Let $f\in\End(\C,0)$ be
tangent to the identity. The element $\mu_f\in\ca M_r$ given by this 
procedure is the {\sl sectorial invariant} of $f$. 

Then the
holomorphic classification proved by
\'Ecalle and Voronin is

\newthm Theorem \EV: (\'Ecalle, 1981 [\'E2--3]; Voronin, 1981 [V]) Let
$f$,~$g\in\End(\C,0)$ be two holomorphic local dynamical systems tangent to the
identity. Then $f$ and $g$ are holomorphically locally conjugated if and only
if they have the same multiplicity, the same index and the same sectorial
invariant. Furthermore, for any $r\ge1$, $\beta\in\C$ and $\mu\in\ca M_r$ there
exists $f\in\End(\C,0)$ tangent to the identity with multiplicity~$r+1$,
index~$\beta$ and sectorial invariant~$\mu$.  

\Rem In particular, holomorphic local dynamical systems tangent to the
identity give examples of local dynamical systems that are topologically
conjugated without being neither holomorphically nor formally conjugated, and of
local dynamical systems that are formally conjugated without being
holomorphically conjugated.

Finally, if $f\in\End(\C,0)$ satisfies $a_1=e^{2\pi i p/q}$, then $f^q$ is
tangent to the identity. Therefore we can apply the previous results to~$f^q$
and then infer informations about the dynamics of the original~$f$,
because of the following

\newthm Lemma \addpar: Let $f$, $g\in\End(\C,0)$ be two holomorphic
local dynamical systems with the same multiplier $e^{2\pi ip/q}\in S^1$.
Then $f$ and $g$ are holomorphically locally conjugated if and only if
$f^q$ and $g^q$ are.

\pf One direction is obvious. For the converse, let $\phe$ be a germ
conjugating $f^q$ and $g^q$; in particular,
$$
g^q=\phe^{-1}\circ f^q\circ\phe=(\phe^{-1}\circ f\circ\phe)^q\;.
$$
So, up to replacing $f$ by $\phe^{-1}\circ f\circ\phe$, we can assume 
that $f^q=g^q$. Put
$$
\psi=\sum_{k=0}^{q-1}g^{q-k}\circ f^k=\sum_{k=1}^q g^{q-k}\circ f^k\;.
$$
The germ $\psi$ is a local biholomorphism, because $\psi'(0)=q\ne 0$,
and it is easy to check that $\psi\circ f=g\circ\psi$.\qedn

We list here a few results; see~[Mi], [Ma], [C], [\'E2--3], [V] and [BH] for
proofs and further details.

\newthm Proposition \tduea: Let $f\in\End(\C,0)$ be a holomorphic local
dynamical system with multiplier $\lambda\in S^1$, and assume that $\lambda$ is a
primitive root of the unity of order $q$. Assume that $f^q\not\equiv\id$. Then
there exist $n\ge 1$ and $\alpha\in\C$ such that $f$ is formally conjugated to
$$
g(z)=\lambda z- z^{nq+1}+ \alpha z^{2nq+1}\;.
$$

\Def The number $n$ is the {\sl parabolic multiplicity} of $f$, and $\alpha\in\C$ is the {\sl index} of~$f$; the {\sl iterative residue} of $f$ is then given by
$$
\hbox{Resit}(f)={nq+1\over 2}-\alpha\;.
$$

\newthm Proposition \ttrea: (Camacho) Let $f\in\End(\C,0)$ be a holomorphic local
dynamical system with multiplier~$\lambda\in S^1$, and assume that $\lambda$ is a
primitive root of the unity of order $q$. Assume that $f^q\not\equiv\id$, and has parabolic multiplicity~$n\ge 1$. Then $f$ is topologically conjugated to
$$
g(z)=\lambda z- z^{nq+1}\;.
$$

\newthm Theorem \tquattroa: (Leau-Fatou) Let $f\in\End(\C,0)$ be a holomorphic 
local dynamical system with multiplier~$\lambda\in S^1$, and assume that $\lambda$ is a
primitive root of the unity of order $q$. Assume that $f^q\not\equiv\id$,
and let $n\ge 1$ be the parabolic multiplicity of~$f$. Then
$f^q$ has multiplicity~$nq+1$, and $f$ acts on
the attracting (respectively, repelling) petals of $f^q$ as a permutation
composed by $n$ disjoint cycles. Finally, $K_f=K_{f^q}$.

Furthermore, it is possible to define the sectorial invariant of such
a holomorphic local dynamical system, composed by $2nq$ germs
whose multipliers still satisfy \eqocon, and the analogue of 
Theorem~\EV\ holds.

\smallsect 4. One complex variable: the elliptic case

We are left with the elliptic case:
$$
f(z)=e^{2\pi i\theta}z+a_2z^2+\cdots\in\C_0\{z\}\;,
\neweq\eqquno
$$
with $\theta\notin\Q$. It turns out that the local dynamics depends mostly
on numerical properties of~$\theta$. The main question here is whether
such a local dynamical system is holomorphically conjugated to its
linear part. Let us introduce a bit of terminology.

\Def We shall say that a holomorphic dynamical system of the form~\eqquno\ is {\sl holomorphically linearizable} if it is
holomorphically locally conjugated to its linear part, the irrational
rotation
$z\mapsto e^{2\pi i\theta}z$. In this case, we shall say that $0$ is a {\sl Siegel point} for~$f$; otherwise, we shall say that it is a {\sl Cremer point.}

It turns out tha for a full measure
subset~$B$ of $\theta\in[0,1]\setminus\Q$ all holomorphic local dynamical systems
of the form~\eqquno\ are holomorphically linearizable. Conversely, the complement $[0,1]\setminus B$ is a
$G_\delta$-dense set, and for all~$\theta\in [0,1]\setminus B$ the quadratic
polynomial~$z\mapsto z^2+e^{2\pi i\theta}z$ is not holomorphically linearizable.
This is the gist of the results due to Cremer, Siegel, Bryuno and Yoccoz we shall describe in this section.

The first worthwhile observation in this setting is that it is possible to give
a topological characterization of holomorphically linearizable local
dynamical systems.

\Def We shall say that $p$
is {\sl stable} for~$f\in\End(M,p)$ if it belongs to the interior of~$K_f$.

\newthm Proposition \local: Let $f\in\End(\C,0)$ be a holomorphic local
dynamical system with multiplier $\lambda\in S^1$. Then
$f$ is holomorphically linearizable if and only if it is topologically
linearizable if and only if $0$ is stable for~$f$.

\pf If $f$ is holomorphically linearizable it is topologically linearizable, and if it
is topologically linearizable (and $|\lambda|= 1$) then it is stable. Assume that
$0$ is stable, and set
$$
\phe_k(z)={1\over k}\sum_{j=0}^{k-1}{f^j(z)\over\lambda^j}\;,
$$
so that $\phe'_k(0)=1$ and 
$$
\phe_k\circ f=\lambda\phe_k+{\lambda\over k}\left({f^k\over\lambda^k}-\id
\right)\;.
\neweq\eqqdue
$$
The stability of~$0$ implies that there are bounded open sets $V\subset U$ containing the origin such that $f^k(V)\subset U$ for all $k\in\N$. Since $|\lambda|=1$, it
follows that $\{\phe_k\}$ is a uniformly bounded family on~$V$, and hence, by Montel's
theorem, it admits a converging subsequence. But
\eqqdue\ implies that a converging subsequence converges to a conjugation between~$f$ and the
rotation~$z\mapsto\lambda z$, and so $f$ is holomorphically linearizable.\qedn

The second important observation is that two elliptic holomorphic local dynamical
systems with the same multiplier are always formally conjugated:

\newthm Proposition \qdue: Let $f\in\End(\C,0)$ be a holomorphic local dynamical
system of multiplier $\lambda=e^{2\pi i\theta}\in S^1$ with $\theta\notin\Q$.
Then $f$ is formally conjugated to its linear part, by a unique formal power
series tangent to the identity.

\pf We shall prove that there is a unique formal power series 
$$
h(z)=z+h_2z^2+\cdots\in\C[[z]]
$$ 
such that $h(\lambda z)=f\bigl(h(z)\bigr)$. Indeed we have
$$
\eqalign{
h(\lambda z)-f\bigl(h(z)\bigr)&=\sum_{j\ge 2}\left\{\bigl[(\lambda^j-\lambda)h_j-a_j\bigr]
z^j-a_j\sum_{\ell=1}^j{j\choose \ell}z^{\ell+j}\left(\sum_{k\ge 2}h_k
z^{k-2}\right)^{\!\ell}\right\}\cr
&=\sum_{j\ge 2}\bigl[(\lambda^j-\lambda)h_j-a_j-X_j(h_2,\ldots,h_{j-1})\bigr]z^j\;,
\cr}
\neweq\eqpass
$$
where $X_j$ is a polynomial in $j-2$ variables with coefficients depending
on~$a_2,\ldots,a_{j-1}$. It follows that the coefficients of~$h$ are
uniquely determined by induction using the formula 
$$
h_j={a_j+X_j(h_2,\ldots,h_{j-1})\over\lambda^j-\lambda}\;.
\neweq\eqqtre
$$
In particular, $h_j$ depends only on $\lambda$,
$a_2,\ldots,a_j$.\qedn

\Rem The same proof shows that any holomorphic local dynamical system with multiplier $\lambda\ne 0$ and not a root of unity is formally conjugated to
its linear part.

The formal power series linearizing~$f$ is not converging if its coefficients
grow too fast. Thus \eqqtre\ links the radius of convergence of~$h$ to the
behavior of~$\lambda^j-\lambda$: if the latter becomes too small, the series
defining~$h$ does not converge. This is known as the {\it small denominators
problem} in this context. 

It is then natural to introduce the following quantity:
$$
\Omega_\lambda(m)=\min_{1\le k\le m}|\lambda^k-\lambda|\;,
$$
for $\lambda\in S^1$ and $m\ge 1$. Clearly, $\lambda$ is a root of unity if and
only if $\Omega_\lambda(m)=0$ for all $m$ greater or equal to some~$m_0\ge 1$;
furthermore, 
$$
\lim_{m\to+\infty}\Omega_\lambda(m)=0
$$ 
for all~$\lambda\in S^1$. 

The first one to actually prove that there are non-linearizable elliptic
holomorphic local dynamical systems has been Cremer, in 1927 [Cr1]. His more
general result is the following:

\newthm Theorem \Cremer: (Cremer, 1938 [Cr2]) Let $\lambda\in S^1$ be
such that
$$
\limsup_{m\to+\infty}{1\over
m}\log{1\over\Omega_\lambda(m)}=+\infty\;.
\neweq\eqqquattro
$$
Then there exists $f\in\End(\C,0)$ with multiplier~$\lambda$ which is not
holomorphically linearizable.
Furthermore, the set of $\lambda\in S^1$ satisfying~$\eqqquattro$ contains
a~$G_\delta$-dense set. 

\pf Choose inductively $a_j\in\{0,1\}$ so
that $|a_j+X_j|\ge1/2$ for all~$j\ge 2$, where~$X_j$ is as in~\eqqtre. Then
$$
f(z)=\lambda z+a_2z^2+\cdots\in\C_0\{z\}\;,
$$
while \eqqquattro\ implies that the radius of convergence of the formal
linearization~$h$ is~0, and thus $f$ cannot be holomorphically linearizable, as
required.

Finally, let $C(q_0)\subset S^1$ denote the set of $\lambda=e^{2\pi i\theta}\in
S^1$ such that
$$
\left|\theta-{p\over q}\right|<{1\over 2^{q!}}
\neweq\eqCC
$$
for some~$p/q\in\Q$ in lowest terms, with $q\ge q_0$. Then it is not difficult to
check that each $C(q_0)$ is a dense open set in~$S^1$, and that all
$\lambda\in\ca C=\bigcap_{q_0\ge 1}C(q_0)$ satisfy~\eqqquattro. Indeed, if
$\lambda=e^{2\pi i\theta}\in\ca C$ we can find $q\in\N$ arbitrarily large such that there is
$p\in\N$ so that \eqCC\ holds. Now, it is easy to see that
$$
|e^{2\pi i t}-1|\le 2\pi|t|
$$
for all $t\in[-1/2,1/2]$. Then let $p_0$ be the integer closest
to~$q\theta$, so that $|q\theta-p_0|\le1/2$. Then we have
$$
|\lambda^q-1|=|e^{2\pi iq\theta}-e^{2\pi ip_0}|=|e^{2\pi
i(q\theta-p_0)}-1|\le2\pi|q\theta-p_0|\le 2\pi|q\theta-p|<{2\pi\over 2^{q!-1}}
$$
for arbitrarily large $q$, and 
\eqqquattro\ follows.
\qedn

On the other hand, Siegel in 1942 gave
a condition on the multiplier ensuring holomorphic linearizability:

\newthm Theorem \Siegel: (Siegel, 1942 [Si]) Let $\lambda\in S^1$ be such
that there exists $\beta> 1$ and $\gamma>0$ so that
$$
\forevery{m\ge 2}{1\over\Omega_\lambda(m)}\le\gamma\, m^\beta.
\neweq\eqqcinque
$$
Then all $f\in\End(\C,0)$ with multiplier~$\lambda$ are
holomorphically linearizable.
Furthermore, the set of $\lambda\in S^1$ satisfying~$\eqqcinque$ for
some~$\beta> 1$ and $\gamma>0$ is of full Lebesgue measure in~$S^1$.

\Rem If $\theta\in[0,1)\setminus\Q$ is algebraic then $\lambda=e^{2\pi i
\theta}$ satisfies \eqqcinque\ for some $\beta>1$ and $\gamma>0$. However,
the set of $\lambda\in S^1$ satisfying \eqqcinque\ is much larger,
being of full measure.

\Rem It is interesting to notice that for generic (in a topological sense)
$\lambda\in S^1$ there is a non-linearizable holomorphic local dynamical system
with multiplier~$\lambda$, while for almost all (in a measure-theoretic sense)
$\lambda\in S^1$ every holomorphic local dynamical system
with multiplier~$\lambda$ is holomorphically linearizable.

Theorem~\Siegel\ suggests the existence of a number-theoretical condition
on~$\lambda$ ensuring that the origin is a Siegel point for any holomorphic
local dynamical system of multiplier~$\lambda$. And indeed this is the content
of the celebrated {\sl Bryuno-Yoccoz theorem:}

\newthm Theorem \BY: (Bryuno, 1965 [Bry1--3], Yoccoz, 1988 [Y1--2]) Let
$\lambda\in S^1$. Then the following statements are equivalent:
{\smallskip 
\item{\rm(i)} the origin is a Siegel point for the quadratic polynomial $f_\lambda(z)=\lambda
z+z^2$; 
\item{\rm(ii)} the origin is a Siegel point for all 
$f\in\End(\C,0)$ with multiplier~$\lambda$;
\item{\rm(iii)} the number $\lambda$
satisfies {\it Bryuno's condition}
$$
\sum_{k=0}^{+\infty}{1\over2^k}\log{1\over\Omega_\lambda(2^{k+1})}<+\infty\;.
\neweq\eqqsei
$$
}

Bryuno, using majorant series as in Siegel's proof of Theorem~\Siegel\ (see also [He] and
references therein) has proved that condition (iii) implies condition (ii). Yoccoz, using a
more geometric approach based on conformal and quasi-conformal geometry, has proved that
(i) is equivalent to (ii), and that (ii) implies (iii), that is that if $\lambda$ does not
satisfy \eqqsei\ then the origin is a Cremer point for some $f\in\End(\C,0)$ with
multiplier~$\lambda$ --- and hence it is a Cremer point for the quadratic polynomial
$f_\lambda(z)$. See also [P9] for related results.

\Rem Condition \eqqsei\ is usually expressed in a different way. Write
$\lambda=e^{2\pi i\theta}$, and let $\{p_k/q_k\}$ be the sequence of rational
numbers converging to~$\theta$ given by the expansion in continued fractions.
Then \eqqsei\ is equivalent to
$$
\sum_{k=0}^{+\infty} {1\over q_k}\log q_{k+1}<+\infty\;,
$$
while \eqqcinque\ is equivalent to
$$
q_{n+1}=O(q_n^\beta)\;,
$$
and \eqqquattro\ is equivalent to
$$
\limsup_{k\to+\infty} {1\over q_k}\log q_{k+1}=+\infty\;.
$$
See [He], [Y2], [Mi] and references therein for details.

\Rem A clear obstruction to the holomorphic linearization of an elliptic
$f\in\End(\C,0)$ with multiplier $\lambda=e^{2\pi i\theta}\in S^2$ is the existence of {\sl small cycles,} that is of periodic
orbits contained in any neighbourhood of the origin. 
Perez-Marco [P1], using Yoccoz's techniques, has shown that when the series
$$
\sum_{k=0}^{+\infty}{\log\log q_{k+1}\over q_k}
$$ 
converges then every germ with multiplier $\lambda$ is either linearizable
or has small cycles, and that when the series diverges there exists 
such germs with a Cremer point but without small cycles.

The complete proof of Theorem~\BY\ is beyond the scope of this survey. We shall limit ourselves
to describe a proof (adapted from [P\"o]) of the implication
(iii)$\Longrightarrow$(ii), to
report two of the easiest results of [Y2], and to illustrate what is the connection
between condition~\eqqsei\ and the radius of convergence of the formal linearizing map.

Let us begin with Bryuno's theorem:

\newthm Theorem \Bry: (Bryuno, 1965 [Bry1--3]) Assume that
$\lambda=e^{2\pi i\theta}\in S^1$ satisfies the Bryuno's condition
$$
\sum_{k=0}^{+\infty}{1\over2^k}\log{1\over\Omega_\lambda(2^{k+1})}<+\infty\;.
\neweq\eqBry
$$
Then the origin is a Siegel point for all 
$f\in\End(\C,0)$ with multiplier~$\lambda$.

\pf
We already know, thanks to Proposition~\qdue, that there exists
a unique formal power series 
$$
h(z)=z+\sum_{k\ge 2} h_kz^k
$$ 
such that
$h^{-1}\circ f\circ h(z)=\lambda z$; we shall prove that $h$
is actually converging. To do so it suffices to show that
$$
\sup_k{1\over |k|} \log|h_k|<\infty\;.
\neweq\eqsei
$$
Since~$f$ is holomorphic in a neighbourhood of the origin, there exists a number~$M>0$ such that~$|a_k|\le M^k$ for~$k\ge 2$; up
to a linear change of coordinates we can assume that $M=1$,
that is $|a_l|\le 1$ for all $k\ge 2$.

Now, $h(\lambda z)=f\bigl(h(z)\bigr)$ yields
$$
\sum_{k\ge 2}(\lambda^k-\lambda)h_k z^k = \sum_{l\ge 2}a_l \left(\sum_{m\ge1}h_m z^m \right)^l\;.
\neweq\eqdp
$$
Therefore
$$
|h_k|\le \eps_k^{-1} \sum_{k_1+\cdots +k_\nu= k \atop \nu\ge 2} |h_{k_1}|\cdots |h_{k_\nu}|\;,
$$
where
$$
\eps_k = |\lambda^k - \lambda|\;.
$$
Define inductively
$$
\alpha_k= \cases{1&if $k=1\;$,\cr
\noalign{\smallskip}
\displaystyle\sum_{k_1+\cdots + k_\nu =k \atop \nu \ge 2} \alpha_{k_1} \cdots \alpha_{k_\nu}& if $k\ge 2\;$,\cr}
$$
and 
$$
\delta_k = \cases{1&if $k=1\;$,\cr
\noalign{\smallskip}
\displaystyle\eps_k^{-1}\max_{k_1+\cdots + k_\nu =k\atop \nu\ge 2} \delta_{k_1}\cdots\delta_{k_\nu},&if $k\ge 2\;$.\cr}
$$
Then it is easy to check by induction that
$$
|h_k|\le \alpha_k\delta_k
$$
for all $k\ge 2$.
Therefore, to establish \eqsei\ it suffices to prove analogous estimates for~$\alpha_k$ and~$\delta_k$.

To estimate~$\alpha_k$, let~$\alpha= \sum_{k \ge 1}\alpha_k t^k$. We have
$$
\alpha - t = \sum_{k\ge 2} \alpha_k t^k= \sum_{k\ge 2}\left(\sum_{j\ge 1}\alpha_j t^j\right)^k= {\alpha^2 \over 1- \alpha}\;.
$$
This equation has a unique holomorphic solution vanishing at zero
$$
\alpha= {t+1 \over 4} \left(1 - \sqrt{1-{8t\over(1+t)^2}}\right)\;,
$$
defined for~$|t|$ small enough. Hence,
$$
\sup_k {1\over k}\log \alpha_k < \infty\;,
$$
as we wanted.

To estimate~$\delta_k$ we have to take care of small divisors.
First of all, for each~$k\ge 2$ we associate to~$\delta_k$ a specific decomposition of the form
$$
\delta_k= \eps_k^{-1}
\delta_{k_1}\cdots \delta_{k_\nu}\;,
\neweq\eqadelta
$$
with~$k>k_1\ge \cdots\ge k_\nu$, $k=k_1 + \cdots + k_\nu$ and $\nu\ge2$,
and hence, by induction, a specific decomposition of the form
$$
\delta_k= \eps_{l_0}^{-1}\eps_{l_1}^{-1}\cdots\eps_{l_q}^{-1}\;,
\neweq\eqdelta
$$
where~$l_0=k$ and $k>l_1\ge\cdots\ge l_q\ge2$. For $m\ge 2$ let 
$N_m(k)$
be the number of factors~$\eps_{l}^{-1}$ in the expression \eqdelta\ 
of~$\delta_k$ satisfying
$$
\eps_{l}<{1\over 4}\,\Omega_\lambda(m)\;.
$$ 
Notice that~$\Omega_\lambda(m)$ is non-increasing with respect to~$m$ and it tends to zero as~$m$ goes to infinity. The next lemma contains the key
estimate.

\newthm Lemma \bruno: For all $m\ge2$ we have
$$
N_m(k)\le \cases{0, & if $k\le m\;$,\cr
\noalign{\smallskip}
\displaystyle {2k\over m}-1,&if $k> m\;$.\cr}
$$

\pf We argue by induction on $k$. If $l\le k\le m$ we have
$\eps_l\ge\Omega_\lambda(m)$, and hence~$N_m(k)=0$.

Assume now~$k>m$, so that~$2k/m -1 \ge 1$. Write $\delta_k$ as in
\eqadelta; we have a few cases to consider.
\smallskip
{\it Case 1:} $\eps_k \ge {1\over4}\,\Omega_\lambda(m)$. Then
$$
N(k) = N(k_1) + \cdots + N(k_\nu)\;,
$$
and applying the induction hypotheses to each term we get~$N(k)\le 
(2k/m) - 1$.
\smallskip
{\it Case 2:} $\eps_k<{1\over4}\,\Omega_\lambda(m)$. Then
$$
N(k) = 1 + N(k_1) + \cdots + N(k_\nu)\;,
$$
and there are three subcases.
\smallskip
{\it Case 2.1:} $k_1\le m$. Then
$$
N(k) = 1 \le {2k\over m} -1\;,
$$
and we are done.
\smallskip
{\it Case 2.2:} $k_1\ge k_2>m$. Then there is~$\nu'$ such that~$2\le\nu'\le\nu$ and~$k_{\nu'}> m\ge k_{\nu'+1}$, and we again have
$$
N(k) = 1 + N(k_1) + \cdots + N(k_{\nu'})\le 1 + {2k\over m} - \nu' \le  
{2k\over m} -1\;.
$$
\smallskip
{\it Case 2.3:} $k_1>m\ge k_2$. Then
$$
N(k)= 1 + N(k_1)\;,
$$
and we have two different subsubcases.
\smallskip
{\it Case 2.3.1:} $k_1\le k-m$. Then
$$
N(k) \le 1 + 2\,{k-m\over m} -1 < {2k\over m} -1\;,
$$
and we are done in this case too.
\smallskip
{\it Case 2.3.2:} $k_1>k-m$. The crucial remark here is that~$\eps_{k_1}^{-1}$ gives no contribute to~$N(k_1)$. Indeed, assume by contradiction
that $\eps_{k_1}<{1\over 4}\,\Omega_\lambda(m)$.
Then
$$
|\lambda^{k_1}|>|\lambda|-{1\over4}\,\Omega_\lambda(m)
\ge 1 - {1\over2} = {1\over2}\;,
$$
because $\Omega_\lambda(m)\le 2$.
Since $k-k_1< m$, it follows that
$$
{1\over2}\,\Omega_\lambda(m)> \eps_k+ \eps_{k_1}
=|\lambda^k -\lambda| + |\lambda^{k_1} -\lambda|
\ge|\lambda^k - \lambda^{k_1}|=|\lambda^{k-k_1} - 1|
\ge \Omega_\lambda(k-k_1+1)\ge\Omega_\lambda(m)\;,
$$
contradiction.

Therefore case~$1$ applies to~$\delta_{k_1}$ and we have
$$
N(k) = 1 + N(k_{1_1}) + \cdots + N(k_{1_{\nu_{1}}})\;,
$$
with~$k> k_1> k_{1_1} \ge \cdots \ge k_{1_{\nu_{1}}}$ and~$k_1 =k_{1_1} + \cdots +k_{1_{\nu_{1}}}$. We can repeat the argument for this decomposition, and we finish unless we run into case~2.3.2 again. However, this loop cannot happen more than~$m+1$ times, and we eventually have to land into a different case. This completes the induction and the proof.\qedn

Let us go back to the proof of Theorem~\Bry. We have to estimate
$$
{1\over k}\log\delta_k = \sum_{j=0}^q {1\over k} \log \eps_{l_j}^{-1}\;.
$$
By Lemma~\bruno, 
$$
\hbox{card}\left\{0\le j\le q\biggm| {1\over 4}\Omega_\lambda(2^{\nu+1}) \le \eps_{l_j} <{1\over4}\,\Omega_\lambda(2^\nu)\right\}
\le N_{2^\nu}(k)\le {2k\over 2^{\nu}} 
$$
for~$\nu\ge 1$. It is also easy to see from the definition of~$\delta_k$ that the number of factors~$\eps_{l_j}^{-1}$ is bounded by~$2k - 1$. In particular,
$$
\hbox{card}\left\{0\le j\le q\biggm| {1\over4}\,\Omega_\lambda(2) \le 
\eps_{l_j} \right\} \le 2k = {2k\over 2^{0}}\;. 
$$
Then
$$
{1\over k} \log \delta_k \le 2 \sum_{\nu\ge 0} {1\over 2^\nu} 
\log\bigl(4\,\Omega_\lambda(2^{\nu +1})^{-1}\bigr) 
	= 2\log 4+2\sum_{\nu \ge 0} {1\over 2^\nu}\log
{1\over\Omega_\lambda(2^{\nu +1})}\;,
$$
and we are done.\qedn

The second result we would like to present is Yoccoz's beautiful
proof of the fact that almost every quadratic polynomial $f_\lambda$
is holomorphically linearizable:

\newthm Proposition \Ypiu: The origin is a Siegel point of $f_\lambda(z)=\lambda z+z^2$ for
almost every $\lambda\in S^1$.

\pf (Yoccoz~[Y2]) The idea is to study the radius of convergence of the inverse
of the linearization  of $f_\lambda(z)=\lambda z+z^2$ when $\lambda\in\Delta^*$. 
Theorem~\Koenigs\ says that there is a unique map $\phe_\lambda$ defined in some
neighbourhood of the origin such that $\phe_\lambda'(0)=1$ and $\phe_\lambda\circ
f=\lambda\phe_\lambda$. Let $\rho_\lambda$ be the radius of convergence
of~$\phe_\lambda^{-1}$; we want to prove that $\phe_\lambda$ is defined in a 
neighbourhood of the unique critical point~$-\lambda/2$ of~$f_\lambda$, and that
$\rho_\lambda=|\phe_\lambda(-\lambda/2)|$.

Let $\Omega_\lambda\subset\subset\C$ be the basin of attraction of the origin, that is the set of $z\in\C$ whose orbit converges to the origin.
Notice that setting $\phe_\lambda(z)=\lambda^{-k}\phe_\lambda\bigl(f_\lambda(z)\bigr)$
we can extend $\phe_\lambda$ to the whole of~$\Omega_\lambda$. Moreover,
since the image of $\phe_\lambda^{-1}$ is contained in~$\Omega_\lambda$, which
is bounded, necessarily $\rho_\lambda<+\infty$. Let
$U_\lambda=\phe_\lambda^{-1}(\Delta_{\rho_\lambda})$. Since we have
$$
(\phe_\lambda'\circ f)f'=\lambda\phe'_\lambda
\neweq\eqYupiu
$$
and $\phe_\lambda$ is invertible in~$U_\lambda$, the function $f$ cannot have critical points
in~$U_\lambda$.

If $z=\phe_\lambda^{-1}(w)\in U_\lambda$, we have $f(z)=\phe^{-1}_\lambda(\lambda w)\in\phe^{-1}
_\lambda(\Delta_{|\lambda|\rho_\lambda})\subset\subset U_\lambda$; therefore 
$$
f(\bar
U_\lambda)\subseteq\bar{f(U_\lambda)}\subset\subset U_\lambda\subseteq\Omega_\lambda,
$$
which implies that $\de U\subset\Omega_\lambda$. So $\phe_\lambda$ is defined on~$\de
U_\lambda$, and clearly $|\phe_\lambda(z)|=\rho_\lambda$ for all $z\in\de U_\lambda$.

If $f$ had no critical points in~$\de
U_\lambda$, \eqYupiu\ would imply that $\phe_\lambda$ has no critical points in~$\de
U_\lambda$. But then $\phe_\lambda$ would be locally invertible in~$\de U_\lambda$, and thus
$\phe^{-1}_\lambda$ would extend across $\de\Delta_{\rho_\lambda}$, impossible.
Therefore $-\lambda/2\in\de U_\lambda$, and $|\phe_\lambda(-\lambda/2)|=\rho_\lambda$, as
claimed.

(Up to here it was classic; let us now start Yoccoz's argument.) Put
$\eta(\lambda)=\phe_\lambda(-\lambda/2)$. From the proof of Theorem~\Koenigs\ one easily sees
that $\phe_\lambda$ depends holomorphically on~$\lambda$; so $\eta\colon\Delta^*\to\C$ is
holomorphic. Furthermore, since $\Omega_\lambda\subseteq\Delta_2$, Schwarz's lemma applied
to $\phe^{-1}_\lambda\colon\Delta_{\rho_\lambda}\to\Delta_2$ yields
$$
1=|(\phe^{-1}_\lambda)'(0)|\le 2/\rho_\lambda,
$$
that is $\rho_\lambda\le 2$. Thus $\eta$ is bounded, and thus it extends holomorphically to
the origin.

So $\eta\colon\Delta\to\Delta_2$ is a bounded holomorphic function not identically zero;
Fatou's theorem on radial limits of bounded holomorphic functions
then implies that 
$$
\rho(\lambda_0):=\limsup_{r\to 1^-}|\eta(r\lambda_0)|>0
$$
for almost every $\lambda_0\in S^1$. This means that we can find $0<\rho_0<\rho(\lambda_0)$
and a sequence $\{\lambda_j\}\subset\Delta$ such that $\lambda_j\to\lambda_0$ and
$|\eta(\lambda_j)|>\rho_0$. This means that $\phe^{-1}_{\lambda_j}$ is defined
in~$\Delta_{\rho_0}$ for all~$j\ge1$; up to a subsequence, we can assume that
$\phe^{-1}_{\lambda_j}\to\psi\colon\Delta_{\rho_0}\to\Delta_2$. But then we have
$\psi'(0)=1$ and 
$$
f_{\lambda_0}\bigl(\psi(z)\bigr)=\psi(\lambda_0 z)
$$
in~$\Delta_{\rho_0}$, and thus the origin is a Siegel point for~$f_{\lambda_0}$.\qedn

The third result we would like to present is the implication (i)${}\Longrightarrow{}$(ii)
in Theorem~\BY. The proof depends on the following result of Douady and Hubbard, obtained
using the theory of quasiconformal maps:

\newthm Theorem \DH: (Douady-Hubbard, 1985 [DH]) Given $\lambda\in\C^*$, let
$f_\lambda(z)=\lambda z+z^2$ be a quadratic polynomial. Then there exists a universal
constant $C>0$ such that for every holomorphic function
$\psi\colon\Delta_{3|\lambda|/2}\to\C$ with
$\psi(0)=\psi'(0)=0$ and $|\psi(z)|\le C|\lambda|$ for all $z\in\Delta_{3|\lambda|/2}$ the
function
$f=f_\lambda+\psi$ is topologically conjugated to~$f_\lambda$
in~$\Delta_{|\lambda|}$.

Then

\newthm Theorem \Ydpiu: (Yoccoz [Y2]) Let $\lambda\in S^1$ be such that the origin is a
Siegel point for $f_\lambda(z)=\lambda z+z^2$. Then the origin is a Siegel point for every
$f\in\End(\C,0)$ with multiplier $\lambda$.

\ifdim\lastskip<\smallskipamount \removelastskip\smallskip\fi
\noindent{\sl Sketch of proof\/}:\enspace Write 
$$
f(z)=\lambda z+a_2 z^2+\sum_{k\ge 3}a_k z^k\;,
$$
and let
$$
f^a(z)=\lambda z+a z^2+\sum_{k\ge 3}a_k z^k\;,
$$
so that $f=f^{a_2}$.
If $|a|$ is large enough then the germ
$$
g^a(z)=af^a(z/a)=\lambda z+z^2+a\sum_{k\ge 3}a_k (z/a)^k=f_\lambda(z)+\psi^a(z)
$$
is defined on $\Delta_{3/2}$ and $|\psi^a(z)|<C$ for all $z\in\Delta_{3/2}$, where $C$ is the constant given by Theorem~\DH. It follows that $g^a$ is topologically
conjugated to~$f_\lambda$. By assumption, $f_\lambda$ is topologically linearizable;
hence $g^a$ is too. Proposition \local\ then implies that $g^a$ is holomorphically linearizable, and hence $f^a$ is too. Furthermore, it is also possible to show
(see, e.g., [BH, Lemma~2.3]) that if $|a|$ is large enough, say $|a|\ge R$, then the
domain of linearization of $g^a$ contains $\Delta_r$, where $r>0$ is such that
$\Delta_{2r}$ is contained in the domain of linearization of~$f_\lambda$. 

So we have proven the assertion if $|a_2|\ge R$; 
assume then $|a_2|<R$.  Since
$\lambda$ is not a root of unity, there exists (Proposition \qdue)
a unique formal power series $\hat
h^a\in\C[[z]]$ tangent to the identity such that $g^a\circ\hat h^a(z)=\hat h^a(\lambda z)$.
If we write
$$
\hat h^a(z)=z+\sum_{k\ge 2} h_k(a)z^k
$$
then $h_k(a)$ is a polynomial in~$a$ of degree~$k-1$, by \eqdp. In particular, 
by the maximum principle we have
$$
|h_k(a_2)|\le \max_{|a|=R}|h_k(a)|
\neweq\eqYtpiu
$$
for all $k\ge 2$. Now, by what we have seen, if $|a|=R$ then $\hat h^a$ is convergent
in a disk of radius $r(a)>0$, and its image contains a disk of radius $r$. Applying
Schwarz's lemma to $(\hat h^a)^{-1}\colon \Delta_r\to\Delta_{r(a)}$ we get $r(a)\ge r$.
But then 
$$
\limsup_{k\to+\infty}|h_k(a_2)|^{1/k}\le \max_{|a|=R}\limsup_{k\to+\infty}|h_k(a)|^{1/k}
={1\over r(a)}\le {1\over r}<+\infty\;;
$$
hence $\hat h^{a_2}$ is convergent, and we are done.\qedn

Finally, we would like to describe the connection between condition \eqqsei\ and linearization. From the
function theoretical side, given $\theta\in[0,1)$ set
$$
r(\theta)=\inf\{r(f)\mid f\in\End(\C,0)\hbox{ has multiplier $e^{2\pi i\theta}$ and it is
defined and injective in $\Delta$}\},
$$
where $r(f)\ge 0$ is the radius of convergence of the unique formal linearization of~$f$
tangent to the identity.

From the number theoretical side, given an irrational number $\theta\in[0,1)$ let $\{p_k/q_k\}$ be
the sequence of rational numbers converging to~$\theta$ given by the expansion in continued
fractions, and put
$$
\eqalign{
\alpha_n&=-{q_n\theta-p_n\over q_{n-1}\theta-p_{n-1}},\qquad \alpha_0=\theta,\cr
\beta_n&=(-1)^n(q_n\theta-p_n),\qquad\beta_{-1}=1.
\cr}
$$

\Def The {\sl Bryuno function} $B\colon[0,1)\setminus\Q\to(0,+\infty]$ is defined by
$$
B(\theta)=\sum_{n=0}^\infty \beta_{n-1}\log{1\over\alpha_n}.
$$

Then Theorem~\BY\ is consequence of what we have seen and the following

\newthm Theorem \Yoccozest: (Yoccoz [Y2]) {\rm (i)} $B(\theta)<+\infty$ if and only if
$\lambda=e^{2\pi i\theta}$ satisfies Bryuno's condition~$\eqqsei$;
\item{\rm(ii)} there exists a universal constant $C>0$ such that
$$
|\log r(\theta)+B(\theta)|\le C
$$
\indent for all $\theta\in[0,1)\setminus\Q$ such that $B(\theta)<+\infty$;
\item{\rm(iii)} if $B(\theta)=+\infty$ then there exists a non-linearizable $f\in\End(\C,0)$
with multiplier~$e^{2\pi i\theta}$.

 If $0$ is a
Siegel point for~$f\in\End(\C,0)$, the local dynamics of~$f$ is completely clear, and simple
enough. On the other hand, if 0 is a Cremer point of~$f$, then the local dynamics of $f$ is
very complicated and not yet completely understood. P\'erez-Marco (in [P2, 4--7]) and
Biswas ([B1, 2]) have studied the topology and the dynamics of the stable set in this case. Some
of their results are summarized in the following

\newthm Theorem \PerezMarco: (P\'erez-Marco, 1995 [P6, 7]) Assume that $0$
is a Cremer point for an elliptic holomorphic local dynamical system
$f\in\End(\C,0)$. Then:
{\smallskip
\item{\rm(i)} The stable set $K_f$ is compact, connected, full (i.e.,
$\C\setminus K_f$ is connected), it is not reduced to~$\{0\}$, and it is not
locally connected at any point distinct from the origin.
\item{\rm(ii)} Any point of~$K_f\setminus\{0\}$ is recurrent (that is, 
a limit point of its orbit).
\item{\rm(iii)} There is an orbit in~$K_f$ which accumulates at the origin, but
no non-trivial orbit converges to the origin.}

\newthm Theorem \biswas: (Biswas, 2007 [B2])
The rotation number and the conformal class of $K_f$ are a complete
set of holomorphic invariants for Cremer points. In other words, two
elliptic non-linearizable holomorphic local dynamical systems $f$ and~$g$
are holomorphically locally conjugated if and only if they have the
same rotation number and there is a biholomorphism of a neighbourhood
of~$K_f$ with a neighbourhood of~$K_g$.

\Rem So, if $\lambda\in S^1$ is not a root of unity and does not
satisfy Bryuno's condition~\eqqsei, we can find $f_1$,~$f_2\in\End(\C,0)$ with
multiplier~$\lambda$ such that~$f_1$ is holomorphically linearizable while~$f_2$
is not. Then $f_1$ and $f_2$ are formally conjugated without being neither
holomorphically nor topologically locally conjugated.

\Rem Yoccoz [Y2] has proved that if $\lambda\in S^1$ is not a root of unity and does 
not satisfy Bryuno's condition~\eqqsei\ then there is an uncountable family of
germs in $\End(\C,O)$ with multiplier~$\lambda$ which are not holomorphically conjugated to each other nor holomorphically conjugated to any entire function.

See also [P1, 3] for other results on the dynamics about a Cremer point.

\smallsect 5. Several complex variables: the hyperbolic case

Now we start the discussion of local dynamics in several complex variables. In
this setting the theory is much less complete than its one-variable counterpart.

\Def Let $f\in\End(\C^n,O)$ be a holomorphic local dynamical system at $O\in\C^n$,
with $n\ge 2$. The {\sl homogeneous expansion} of $f$ is 
$$
f(z)=P_1(z)+P_2(z)+\cdots\in\C_0\{z_1,\ldots,z_n\}^n\;,
$$
where $P_j$ is an $n$-uple of homogeneous polynomials of degree~$j$. In
particular,~$P_1$ is the differential~$df_O$ of~$f$ at the origin, and $f$ is
locally invertible if and only if~$P_1$ is invertible.

We have seen that in dimension one the multiplier (i.e., the derivative at the
origin) plays a main r\^ole. When $n>1$, a similar r\^ole is played by the
eigenvalues of the differential. 

\Def Let $f\in\End(\C^n,O)$ be a holomorphic local dynamical system at $O\in\C^n$,
with $n\ge 2$. Then:
\smallskip
\item{--} if all eigenvalues of $df_O$ have modulus less than 1, we say that the
fixed point~$O$ is {\sl attracting;}
\item{--} if all eigenvalues of $df_O$ have modulus greater than 1, we say that
the fixed point~$O$ is {\sl repelling;}
\item{--} if all eigenvalues of $df_O$ have modulus different from 1, we say that
the fixed point~$O$ is {\sl hyperbolic} (notice that we allow the eigenvalue
zero);
\item{--} if $O$ is attracting or repelling, and $df_O$ is invertible, we say that $f$ is in the {\sl Poincar\'e domain;}
\item{--} if $O$ is hyperbolic, $df_O$ is invertible, and $f$ is not in
the Poincar\'e domain (and thus $df_O$ has both eigenvalues inside the
unit disk and outside the unit disk) we say that $f$ is in the {\sl Siegel
domain;}
\item{--} if all eigenvalues of $df_O$ are roots of unity, we say that the
fixed point~$O$ is {\sl parabolic;} in particular, if $df_O=\id$ we say that $f$
is {\sl tangent to the identity;}
\item{--} if all eigenvalues of $df_O$ have modulus 1 but none is a root of
unity, we say that the fixed point~$O$ is {\sl elliptic;}
\item{--} if $df_O=O$, we say that the fixed point~$O$ is {\sl superattracting.}

\noindent Other cases are clearly possible, but for our aims this list is
enough. In this survey we shall be mainly concerned with hyperbolic and
parabolic fixed points; however, in the last section we shall also present some
results valid in other cases.

Let us begin assuming that the origin is a hyperbolic
fixed point for an $f\in\End(\C^n,O)$ not necessarily invertible. We then have a
canonical splitting
$$
\C^n=E^s\oplus E^u\;,
$$
where $E^s$ (respectively, $E^u$) is the direct sum of the generalized
eigenspaces associated to the eigenvalues of~$df_O$ with modulus less
(respectively, greater) than 1. Then the first main result in this subject is
the famous {\sl stable manifold theorem} (originally due to Perron~[Pe] and
Hadamard~[H]; see [FHY, HK, HPS, Pes, Sh, AM] for proofs in
the $C^\infty$ category, Wu~[Wu] for a proof in the holomorphic category, and
[A3] for a proof in the non-invertible case):

\newthm Theorem \stable: Let $f\in\End(\C^n,O)$ be a holomorphic local dynamical
system with a hyperbolic fixed point at the origin. Then:
{\smallskip
\item{\rm(i)}the stable set~$K_f$ is an embedded complex submanifold of (a
neighbourhood of the origin in)~$\C^n$, tangent to~$E^s$ at the origin;
\item{\rm(ii)}there is an embedded complex submanifold~$W_f$ of (a
neighbourhood of the origin in)~$\C^n$, called the\break\indent {\it unstable
set} of~$f$, tangent to~$E^u$ at the origin, such that $f|_{W_f}$ is
invertible, $f^{-1}(W_f)\subseteq W_f$, and $z\in W_f$\break\indent if and only
if there is a sequence $\{z_{-k}\}_{k\in\N}$ in the domain of~$f$ such that
$z_0=z$ and $f(z_{-k})=z_{-k+1}$ for\break\indent all~$k\ge 1$. Furthermore, if
$f$ is invertible then $W_f$ is the stable set of~$f^{-1}$. }

The proof is too involved to be summarized here; it suffices to say that both
$K_f$ and $W_f$ can be recovered, for instance, as fixed points of a suitable
contracting operator in an infinite dimensional space (see the references quoted
above for details).

\Rem If the origin is an attracting fixed point, then $E^s=\C^n$, and $K_f$ is an
open neighbourhood of the origin, its {\sl basin of attraction.} However, as we
shall discuss below, this does not imply that
$f$ is holomorphically linearizable, not even when it is invertible.
Conversely, if the origin is a repelling fixed point, then
$E^u=\C^n$, and~$K_f=\{O\}$. Again, not all holomorphic local dynamical systems
with a repelling fixed point are holomorphically linearizable.

If a point in the domain~$U$ of a holomorphic local dynamical system with a
hyperbolic fixed point does not belong either to the stable set or to the
unstable set, it escapes both in forward time (that is, its orbit escapes) and in
backward time (that is, it is not the end point of an infinite orbit contained
in~$U$). In some sense, we can think of the stable and unstable sets (or, as
they are usually called in this setting, stable and unstable {\it manifolds\/})
as skewed coordinate planes at the origin, and the orbits outside these
coordinate planes follow some sort of hyperbolic path, entering and leaving any
neighbourhood of the origin in finite time.

Actually, this idea of straightening stable and unstable manifolds can be
brought to fruition (at least in the invertible case), and it yields one of the
possible proofs (see~[HK, Sh, A3] and references therein) of the
{\sl Grobman-Hartman theorem:}

\newthm Theorem \GrobmanHartman: (Grobman, 1959 [G1--2]; Hartman, 1960 [Har])
Let
$f\in\End(\C^n,O)$ be a locally invertible holomorphic local dynamical system
with a hyperbolic fixed point. Then $f$ is topologically locally conjugated to
its differential~$df_O$.

Thus, at least from a topological point of view, the local dynamics about an
invertible hyperbolic fixed point is completely clear. This is definitely not the
case if the local dynamical system is not invertible in a neighbourhood of the
fixed point. For instance, already Hubbard and Papadopol [HP] noticed that
a B\"ottcher-type theorem for superattracting points in several complex
variables is just not true: there are holomorphic local dynamical systems with a
superattracting fixed point which are not even topologically locally conjugated
to the first non-vanishing term of their homogeneous expansion. Recently, Favre
and Jonsson (see, e.g., [Fa] and [FJ1, 2]) have begun a very detailed study of superattracting fixed
points in~$\C^2$, study that might lead to their topological classification.
We shall limit ourselves to quote one result.

\Def Given $f\in\End(\C^2,O)$, we shall denote by $\hbox{\rm Crit}(f)$ the
set of critical points of~$f$. Put
$$
\hbox{\rm Crit}^\infty(f)=\bigcup_{k\ge 0}f^{-k}\bigl(\hbox{\rm Crit}(f)
\bigr)\;;
$$
we shall say that $f$ is {\sl rigid} if (as germ in the origin) 
$\hbox{\rm Crit}^\infty(f)$ is either empty, a smooth curve, or the union of
two smooth curves crossing transversally at the origin. Finally, we
shall say that $f$ is {\sl dominant} if~$\det(df)\not\equiv 0$.

Rigid germs have been classified by Favre [Fa], which isthe reason why
next theorem can be useful for classifying superattracting dynamical
systems: 
 
\newthm Theorem \FJ: (Favre-Jonsson, 2007 [FJ2]) Let $f\in\End(\C^2,O)$
be superattracting and dominant. Then there exist:
\smallskip
\item{\rm (a)} a 2-dimensional
complex manifold~$M$ (obtained by blowing-up a finite number of points; 
see next section);
\item{\rm (b)} a surjective holomorphic map
$\pi\colon M\to\C^2$ such that the restriction $\pi|_{M\setminus E}
\colon M\setminus E\to\C^2\setminus\{O\}$ is a biholomorphism,
where $E=\pi^{-1}(O)$;
\item{\rm (c)} a point $p\in E$; and
\item{\rm(d)} a rigid holomorphic germ $\tilde f\in\End(M,p)$
\smallskip
\noindent so that $\pi\circ\tilde f=f\circ\pi$.

Coming back to hyperbolic dynamical systems,
the holomorphic and even the formal classification are not as simple as the
topological one. The main problem is caused by resonances.

\Def Let $f\in\End(\C^n,O)$ be a holomorphic local dynamical system, and let
denote by
$\lambda_1,\ldots,\lambda_n\in\C$ the eigenvalues of~$df_O$. A 
{\sl resonance} for $f$ is a relation of the form 
$$
\lambda_1^{k_1}\cdots\lambda_n^{k_n}-\lambda_j=0
\neweq\eqcuno
$$
for some $1\le j\le n$ and some $k_1,\ldots,k_n\in\N$ with $k_1+\cdots+k_n\ge
2$. When $n=1$ there
is a resonance if and only if the multiplier is a root of unity, or
zero; but if
$n>1$ resonances may occur in the hyperbolic case too. 

Resonances are the obstruction to formal linearization. Indeed, 
a computation
completely analogous to the one yielding Proposition~\qdue\ shows that 
the coefficients of a
formal linearization have in the denominators quantities of the form
$\lambda_1^{k_1}\cdots\lambda_n^{k_n}-\lambda_j$; hence

\newthm Proposition \qtre: Let $f\in\End(\C^n,O)$ be a locally invertible
holomorphic local dynamical system with a hyperbolic fixed point and no
resonances. Then $f$ is formally conjugated to its differential~$df_O$.

In presence of resonances, even the formal classification is not that easy.
Let us assume, for simplicity, that $df_O$ is in Jordan form, that is
$$
P_1(z)=(\lambda_1z,\epsilon_2 z_1+\lambda_2z_2,\ldots,\epsilon_n
z_{n-1}+\lambda_nz_n)\;,
$$
with $\epsilon_1,\ldots,\epsilon_{n-1}\in\{0,1\}$. 

\Def We shall say that a
monomial
$z_1^{k_1}\cdots z_n^{k_n}$ in the $j$-th coordinate of~$f$ is {\sl resonant} if
$k_1+\cdots+k_n\ge 2$ and
$\lambda_1^{k_1}\cdots\lambda_n^{k_n}=\lambda_j$. 

Then Proposition~\qtre\ can be generalized to (see [Ar, p.~194] or [IY, p.~53] for a proof): 

\newthm Proposition \PoincareDulac: (Poincar\'e [Po], Dulac [Du]) Let $f\in\End(\C^n,O)$ be a locally 
invertible holomorphic local dynamical system with a hyperbolic fixed point.
Then it is formally conjugated to a $g\in\C_0[[z_1,\ldots,z_n]]^n$ such that
$dg_O$ is in Jordan normal form, and $g$ has only resonant monomials.

\Def The formal series~$g$ is called a {\sl Poincar\'e-Dulac normal form} of~$f$. 

The problem with Poincar\'e-Dulac normal forms is that they are not unique. In
particular, one may wonder whether it could be possible to have such a normal
form including {\it finitely many} resonant monomials only (as happened, for
instance, in Proposition~\formaltangent). 

This is
indeed the case (see, e.g., Reich~[Re1]) when $f$ belongs to the 
Poincar\'e domain, that is when $df_O$ is invertible and $O$ is either
attracting or repelling. As far as I know, the
problem of finding canonical formal normal forms 
when
$f$ belongs to the Siegel domain is still open.

It should be remarked that, in the hyperbolic case, the problem of formal
linearization is equivalent to the problem of smooth linearization. This has
been proved by Sternberg~[St1--2] and Chaperon~[Ch]:

\newthm Theorem \Sternberg: (Sternberg, 1957 [St1--2]; Chaperon, 1986 [Ch])
Let 
$f$,~$g\in\End(\C^n,O)$ be two holomorphic local dynamical
systems, and assume that $f$ is locally invertible and with a hyperbolic fixed
point at the origin. Then $f$ and $g$ are formally conjugated if and only if
they are smoothly locally conjugated. In particular, $f$ is smoothly
linearizable if and only if it is formally linearizable. Thus if there are
no resonances then $f$ is smoothly linearizable.

Even without resonances, the holomorphic linearizability is not guaranteed.
The easiest positive result is due to Poincar\'e~[Po] who, using majorant
series, proved the following

\newthm Theorem \Poincare: (Poincar\'e, 1893 [Po])  Let $f\in\End(\C^n,O)$ be
a  locally invertible holomorphic local dynamical system in the Poincar\'e 
domain. Then $f$ is holomorphically linearizable if and
only if it is formally linearizable. In particular, if there are no resonances
then $f$ is holomorphically linearizable.

Reich~[Re2] describes holomorphic normal forms when $df_O$ belongs to the
Poincar\'e domain and there are resonances (see also~[\'EV]); P\'erez-Marco [P8] 
discusses the problem of holomorphic linearization in the presence of resonances.

When $df_O$ belongs to the Siegel domain, even without resonances, the formal
linearization might diverge. To describe the known results, let us introduce the
following definition:

\Def For $\lambda_1,\ldots,\lambda_n\in\C$ and $m\ge 2$ set
$$
\Omega_{\lambda_1,\ldots,\lambda_n}(m)=\min\bigl\{|\lambda_1^{k_1}\cdots
\lambda_n^{k_n}-\lambda_j|\bigm| k_1,\ldots,k_n\in\N,\,2\le k_1+\cdots+k_n\le m, 
\,1\le j\le n\bigr\}\;.
\neweq\eqcB
$$
If $\lambda_1,\ldots,\lambda_n$ are the eigenvalues of~$df_O$, we shall
write~$\Omega_f(m)$ for~$\Omega_{\lambda_1,\ldots,\lambda_n}(m)$.

It is clear that $\Omega_f(m)\ne 0$ for all~$m\ge 2$ if and only if there are no
resonances. It is also not difficult to prove that if $f$ belongs to the
Siegel domain then
$$
\lim_{m\to+\infty}\Omega_f(m)=0\;,
$$
which is the reason why, even without resonances, the formal linearization might
be diverging, exactly as in the one-dimensional case. As far as I know, the best
positive and negative results in this setting are due to Bryuno~[Bry2--3] 
(see also [R\"u]), and
are a natural generalization of their one-dimensional counterparts, whose proofs are
obtained adapting the proofs of Theorems~\Bry\ and \Cremer\ respectively:

\newthm Theorem \Bryunohyp: (Bryuno, 1971 [Bry2--3]) Let $f\in\End(\C^n,O)$ be a
holomorphic local dynamical system such that
$f$ belongs to the Siegel domain, has no resonances, and $df_O$ is diagonalizable.
Assume moreover that 
$$
\sum_{k=0}^{+\infty}{1\over 2^k}\log{1\over\Omega_f(2^{k+1})}<+\infty\;.
\neweq\eqcdue
$$
Then $f$ is holomorphically linearizable.

\newthm Theorem \BryunoCremer: Let $\lambda_1,\ldots,\lambda_n\in\C$ be without
resonances and such that
$$
\limsup_{m\to+\infty}{1\over
m}\log{1\over\Omega_{\lambda_1,\ldots,\lambda_n}(m)}=+\infty\;.
$$
Then there exists $f\in\End(\C^n,O)$, with $df_O=\hbox{\rm
diag}(\lambda_1,\ldots,\lambda_n)$, not holomorphically 
linearizable.

\Rem These theorems hold even without hyperbolicity assumptions. 

\Rem It should be remarked that, contrarily to the one-dimensional case, it is not yet
known whether condition \eqcdue\ is necessary for the holomorphic 
linearizability of all holomorphic local dynamical systems with a given linear
part belonging to the Siegel domain. However, it is easy to check that if $\lambda\in S^1$
does not satisfy the one-dimensional Bryuno condition then any $f\in\End(\C^n,O)$ of
the form
$$
f(z)=\bigl(\lambda z_1+z_1^2, g(z)\bigr)
$$
is not holomorphically linearizable: indeed, if $\phe\in\End(\C^n,O)$ is a holomorphic linearization of~$f$, then $\psi(\zeta)=\phe(\zeta,O)$ is a holomorphic linearization
of the quadratic polynomial $\lambda z+z^2$, against Theorem~\BY.

P\"oschel~[P\"o] shows how to modify \eqcB\ and \eqcdue\ to get partial linearization results
along submanifolds, and [R] (see also [Ro1]) explains when it is possible to pass from a partial linearization
to a complete linearization even in presence of resonances.
See also Russmann [R\"u] for another proof of Theorem~\Bryunohyp, and
Il'yachenko~[I1] for an important
result related to Theorem~\BryunoCremer. Finally, in [DG] are discussed results
in the spirit of Theorem~\Bryunohyp\ without assuming that the differential is
diagonalizable.

\smallsect 6. Several complex variables: the parabolic case

A first natural question in the several complex variables parabolic case is
whether a result like the Leau-Fatou flower theorem holds, and, if so, in
which form. To present what is known on this subject in this section we shall
restrict our attention to holomorphic local dynamical systems 
tangent to the identity; consequences on dynamical systems with a more general
parabolic fixed point can be deduced taking a suitable iterate (but see also
the end of this section for results valid when the differential at the
fixed point is not diagonalizable).

So we are interested in the local dynamics of a holomorphic local dynamical
system $f\in\End(\C^n,O)$ of the form
$$
f(z)=z+P_\nu(z)+P_{\nu+1}(z)+\cdots\in\C_0\{z_1,\ldots,z_n\}^n\;,
\neweq\eqpsc
$$
where $P_\nu$ is the first non-zero term in the homogeneous expansion of~$f$.

\Def If $f\in\End(\C^n,O)$ is of the form \eqpsc,
the number~$\nu\ge 2$ is the {\sl order} of~$f$. 

The two main ingredients in the statement of the Leau-Fatou flower theorem were
the attracting directions and the petals. Let us first describe
a several variables analogue of attracting directions.

\Def Let $f\in\End(\C^n,O)$ be tangent at the identity and of order~$\nu$. A {\sl
characteristic direction} for~$f$ is a non-zero vector $v\in\C^n\setminus\{O\}$
such that $P_\nu(v)=\lambda v$ for some~$\lambda\in\C$. If $P_\nu(v)=O$ (that
is, $\lambda=0$) we shall say that $v$ is a {\sl degenerate} characteristic
direction; otherwise, (that is, if $\lambda\ne 0$) we shall say that $v$ is {\sl
non-degenerate.} We shall say that $f$ is {\sl dicritical} if all directions are characteristic; {\sl non-dicritical} otherwise.

\Rem It is easy to check that $f\in\End(\C^n,O)$ of the form \eqpsc\ is dicritical if
and only if $P_\nu\equiv\lambda\id$, where $\lambda\colon\C^n\to\C$ is
a homogeneous polynomial of degree $\nu-1$. In particular, generic germs 
tangent to the identity are non-dicritical.   

\Rem There is an equivalent definition of characteristic directions that shall be
useful later on. The $n$-uple of $\nu$-homogeneous polynomials~$P_\nu$ induces a
meromorphic self-map of~$\P^{n-1}(\C)$, still denoted by~$P_\nu$. Then, under
the canonical projection $\C^n\setminus\{O\}\to\P^{n-1}(\C)$ non-degenerate characteristic directions correspond
exactly to fixed points of~$P_\nu$, and degenerate characteristic directions
correspond exactly to indeterminacy points of~$P_\nu$. In generic cases, there is only
a finite number of characteristic directions, and using
Bezout's theorem it is easy to prove (see, e.g., [AT1]) that this number, counting according to a suitable multiplicity,
is given by~$(\nu^n-1)/(\nu-1)$.

\Rem The characteristic directions are {\it complex} directions; in particular,
it is easy to check that $f$ and $f^{-1}$ have the same characteristic
directions. Later on we shall see how to associate to (most) characteristic
directions $\nu-1$ petals, each one in some sense centered about a {\it
real} attracting direction corresponding to the same complex characteristic
direction.

The notion of characteristic directions has a dynamical origin. 

\Def We shall say
that an orbit
$\{f^k(z_0)\}$ converges to the origin {\sl tangentially} to a
direction~$[v]\in\P^{n-1}(\C)$ if $f^k(z_0)\to O$ in~$\C^n$ and
$[f^k(z_0)]\to[v]$ in~$\P^{n-1}(\C)$, where $[\cdot]\colon\C^n
\setminus\{O\}\to\P^{n-1}(\C)$ denotes the canonical projection.

Then

\newthm Proposition \suno: Let $f\in\End(\C^n,O)$ be a holomorphic dynamical
system tangent to the identity. If there is an orbit of~$f$ converging to the
origin tangentially to a direction~$[v]\in\P^{n-1}(\C)$, then $v$ is a
characteristic direction of~$f$.

\ifdim\lastskip<\smallskipamount \removelastskip\smallskip\fi
\noindent{\sl Sketch of proof\/}:\enspace ([Ha2]) For simplicity let us
assume~$\nu=2$; a similar argument works for~$\nu>2$.  

If $v$ is a degenerate characteristic direction, there is nothing to prove. If
not, up to a linear change of coordinates we can assume $[v]=[1:v']$ and write
$$
\cases{f_1(z)=z_1+p^1_2(z_1,z')+p_3^1(z_1,z')+\cdots\;,\cr
\noalign{\smallskip}
f'(z)=z'+p'_2(z_1,z')+p'_3(z_1,z')+\cdots\;,\cr}
$$
where $z'=(z_2,\ldots,z_n)\in\C^{n-1}$, $f=(f_1,f')$, $P_j=(p_j^1,p'_j)$ and so
on, with $p_2^1(1,v')\ne 0$. Making the substitution 
$$
\cases{w_1=z_1\;,\cr
z'=w'z_1\;,\cr}
\neweq\eqblowup
$$
which is a change of variable off the hyperplane $z_1=0$, the map~$f$ becomes
$$
\cases{\tilde f_1(w)=w_1+p^1_2(1,w')w_1^2+p_3^1(1,w')w_1^3+\cdots\;,\cr
\noalign{\smallskip}
\tilde f'(w)=w'+r(w')w_1+O(w_1^2)\;,\cr}
\neweq\eqsdue
$$
where $r(w')$ is a polynomial such that $r(v')=O$ if and only if $[1:v']$ is a
characteristic direction of~$f$ with~$p_2^1(1,v')\ne0$.

Now, the hypothesis is that there exists an orbit $\{f^k(z_0)\}$ converging to
the origin and such that $[f^k(z_0)]\to[v]$. Writing $\tilde
f^k(w_0)=\bigl(w_1^k,(w')^k\bigr)$, this implies that $w_1^k\to 0$ and
$(w')^k\to v'$. Then, arguing as in the proof of \eqtquattroa, it is not difficult to prove that 
$$
\lim_{k\to+\infty}{1\over kw_1^k}=-p_2^1(1,v')\;,
$$
and then that $(w')^{k+1}-(w')^k$ is of order $r(v')/k$. This implies
$r(v')=O$, as claimed, because otherwise the telescopic series 
$$
\sum_k\bigl((w')^{k+1}-(w')^k\bigr)
$$
would not be convergent.\qedn

\Rem There are examples of germs $f\in\End(\C^2,O)$ tangent to
the identity with orbits converging to the origin without being tangent to any
direction: for instance 
$$
f(z,w)=\bigl(z+\alpha zw,w+\beta w^2+o(w^2)\bigr)
$$ 
with $\alpha$, $\beta\in\C^*$, $\alpha\ne\beta$ and
$\Re(\alpha/\beta)=1$ (see [Ri1] and [AT3]).

The several variables analogue of a petal is given by the
notion of parabolic curve. 

\Def A {\sl parabolic curve} for $f\in\End(\C^n,O)$
tangent to the identity is an injective holomorphic map
$\phe\colon\Delta\to\C^n\setminus\{O\}$ satisfying the following properties:
\smallskip
\item{(a)} $\Delta$ is a simply connected domain in~$\C$ with $0\in\de\Delta$;
\item{(b)} $\phe$ is continuous at the origin, and $\phe(0)=O$;
\item{(c)} $\phe(\Delta)$ is $f$-invariant, and $(f|_{\phe(\Delta)})^k\to O$
uniformly on compact subsets as $k\to+\infty$.
\smallskip
\noindent Furthermore, if $[\phe(\zeta)]\to[v]$ in $\P^{n-1}(\C)$ as~$\zeta\to
0$ in~$\Delta$, we shall say that the parabolic curve $\phe$ is {\sl tangent} to
the direction~$[v]\in\P^{n-1}(\C)$. 

Then the first main generalization of the Leau-Fatou flower theorem to several
complex variables is due to \'Ecalle and Hakim (see also Weickert [W]):

\newthm Theorem \EcalleHakim: (\'Ecalle, 1985 [\'E4]; Hakim, 1998 [Ha2]) Let
$f\in\End(\C^n,O)$ be a holomorphic local dynamical system tangent to the
identity of order~$\nu\ge 2$. Then for any non-degenerate characteristic
direction~$[v]\in\P^{n-1}(\C)$ there exist (at least) $\nu-1$ parabolic curves
for~$f$ tangent to~$[v]$.

\ifdim\lastskip<\smallskipamount \removelastskip\smallskip\fi
\noindent{\sl Sketch of proof\/}:\enspace \'Ecalle proof is based on his theory
of resurgence of divergent series; we shall describe here the ideas behind
Hakim's proof, which depends on more standard arguments. 

For the sake of simplicity, let us assume $n=2$; without loss of generality we
can also assume $[v]=[1:0]$. Then after a linear change of variables and
a transformation of the kind~\eqblowup\ we are reduced
to prove the existence of a parabolic curve at the origin for a map of the form
$$
\cases{f_1(z)=z_1-z_1^\nu+O(z_1^{\nu+1},z_1^\nu z_2)\;,\cr
\noalign{\smallskip}
f_2(z)=z_2\bigl(1-\lambda
z_1^{\nu-1}+O(z_1^\nu,z_1^{\nu-1} z_2)\bigr)+z_1^\nu
	\psi(z)\;,\cr}
$$
where $\psi$ is holomorphic with $\psi(O)=0$, and $\lambda\in\C$.
Given
$\delta>0$, set
$D_{\delta,\nu}=\{\zeta\in\C\mid |\zeta^{\nu-1}-\delta|<\delta\}$; this open set
has $\nu-1$ connected components, all of them satisfying condition (a) on the
domain of a parabolic curve. Furthermore, if $u$ is a holomorphic function
defined on one of these connected components, of the form $u(\zeta)=\zeta^2
u_o(\zeta)$ for some bounded holomorphic function~$u_o$, and such that
$$
u\bigl(f_1\bigl(\zeta,u(\zeta)\bigr)\bigr)=f_2\bigl(\zeta,u(\zeta)\bigr)\;,
\neweq\eqstre
$$
then it is not difficult to verify that $\phe(\zeta)=\bigl(\zeta,u(\zeta)\bigl)$
is a parabolic curve for~$f$ tangent to~$[v]$. 

So we are reduced to finding a solution of \eqstre\ in each connected component
of~$D_{\delta,\nu}$, with $\delta$ small enough. For any holomorphic $u=\zeta^2
u_o$ defined in such a connected component, let
$f_u(\zeta)=f_1\bigl(\zeta,u(\zeta)\bigr)$, put
$$
H(z)=z_2-{z_1^\lambda\over f_1(z)^\lambda} f_2(z)\;,
$$
and define the operator $T$ by setting
$$
(Tu)(\zeta)=\zeta^\lambda \sum_{k=0}^\infty {H\bigl(f_u^k(\zeta), u\bigl(
f_u^k(\zeta)\bigr)\bigr)\over f_u^k(\zeta)^\lambda}\;.
$$
Then, if $\delta>0$ is small enough, it is possible to prove that $T$ is
well-defined, that $u$ is a fixed point of~$T$ if and only if it satisfies
\eqstre, and that $T$ is a contraction of a closed convex set of a
suitable complex Banach space --- and thus it has a fixed point. In this way if
$\delta>0$ is small enough we get a unique solution of \eqstre\ for each
connected component of~$D_{\delta,\nu}$, and hence $\nu-1$ parabolic curves
tangent to~$[v]$.\qedn

\Def A set of $\nu-1$ parabolic curves obtained in this way is a {\sl
Fatou flower} for~$f$ tangent to~$[v]$. 

\Rem When there is a one-dimensional $f$-invariant complex submanifold passing
through the origin tangent to a characteristic direction~$[v]$, the
previous theorem is just a consequence of the usual one-dimensional theory. But
it turns out that in most cases such an $f$-invariant complex submanifold does
not exist: see [Ha2] for a concrete example, and [\'E4] for a general
discussion. 

We can also have $f$-invariant complex submanifolds of
dimension strictly greater than one attracted by the origin. 

\Def Given
a holomorphic local dynamical system $f\in\End(\C^n,O)$ tangent to the
identity of order~$\nu\ge 2$, and
a non-degenerate characteristic direction $[v]\in\P^{n-1}(\C)$, 
the eigenvalues $\alpha_1,\ldots,\alpha_{n-1}\in\C$ of the linear operator
$d(P_\nu)_{[v]}-\id\colon T_{[v]}\P^{n-1}(\C)\to T_{[v]}\P^{n-1}(\C)$
are the {\sl directors} of~$[v]$. 

Then, using a more elaborate
version of her proof of Theorem~\EcalleHakim, Hakim has been able to prove the
following:

\newthm Theorem \Hakim: (Hakim, 1997 [Ha3]) Let $f\in\End(\C^n,O)$ be a
holomorphic local dynamical system tangent to the identity of order~$\nu\ge 2$.
Let
$[v]\in\P^{n-1}(\C)$ be a non-degenerate characteristic direction, with directors
$\alpha_1,\ldots,\alpha_{n-1}\in\C$.
Furthermore, assume that $\Re\alpha_1,\ldots,\Re\alpha_d>0$ and
$\Re\alpha_{d+1},\ldots,\Re\alpha_{n-1}\le 0$ for a suitable $d\ge 0$. Then:
\smallskip
{\item{\rm(i)}There exists an $f$-invariant $(d+1)$-dimensional complex
submanifold~$M$ of~$\C^n$, with the origin in its boundary, such that the orbit
of every point of~$M$ converges to the origin tangentially to~$[v]$;
\item{\rm(ii)} $f|_M$ is holomorphically conjugated to the translation
$\tau(w_0,w_1,\ldots,w_d)=(w_0+1,w_1,\ldots,w_d)$ defined\break\indent on a
suitable right half-space in~$\C^{d+1}$. }

\Rem In particular, if all the directors of~$[v]$
have positive real part, there is an open domain attracted by the origin.
However, the condition given by Theorem~\Hakim\ is not necessary for the
existence of such an open domain; see Rivi~[Ri1] for an easy example, and
Ushiki~[Us] for a more elaborate example with an open domain attracted by the
origin where $f$ cannot be conjugate to a translation.

In his monumental work [\'E4] \'Ecalle has given a complete set of formal
invariants for holomorphic local dynamical systems tangent to the identity with
at least one non-degenerate characteristic direction. For instance, he
has proved the following

\newthm Theorem \Ecalleformal: (\'Ecalle, 1985 [\'E4]) Let $f\in\End(\C^n,O)$
be a holomorphic local dynamical system tangent to the identity of order~$\nu\ge
2$. Assume that
{\smallskip
\item{\rm(a)} $f$ has exactly $(\nu^n-1)/(\nu-1)$ distinct non-degenerate
characteristic directions and no degenerate characteristic directions;
\item{\rm(b)} the directors of any non-degenerate characteristic
direction are irrational and mutually independent over~$\Z$.
\smallskip
\noindent Choose a non-degenerate characteristic direction~$[v]\in\P^{n-1}(\C)$,
and let $\alpha_1,\ldots,\alpha_{n-1}\in\C$ be its directors. Then there exist a
unique $\rho\in\C$ and unique (up to dilations) formal
series~$R_1,\ldots,R_n\in\C[[z_1,\ldots,z_n]]$, where each $R_j$ contains only
monomial of total degree at least~$\nu+1$ and of partial degree in~$z_j$ at
most~$\nu-2$, such that $f$ is formally conjugated to the time-1 map of the
formal vector field
$$
X={1\over(\nu-1)(1+\rho z_n^{\nu-1})}\left\{[-z_n^\nu+R_n(z)]{\de\over\de
z_n}+\sum_{j=1}^{n-1}[-\alpha_j z_n^{\nu-1}z_j+R_j(z)]{\de\over\de z_j}\right\}.
$$}

Other approaches to the formal classification, at least in dimension~2, are
described in~[BM] and in~[AT2].

Using his theory of resurgence, and always assuming the
existence of at least one non-degenerate characteristic direction, \'Ecalle has
also provided a set of holomorphic invariants for holomorphic local dynamical
systems tangent to the identity, in terms of differential operators with formal
power series as coefficients. Moreover, if the directors of all
non-degenerate characteristic directions are irrational and satisfy a suitable
diophantine condition, then these invariants become a complete set of
invariants. See [\'E5] for a description of his results, and [\'E4]
for the details.

Now, all these results beg the question: what happens when there are no
non-degenerate characteristic directions? For instance, this is the case for
$$
\cases{f_1(z)=z_1+bz_1z_2+z_2^2,\cr
\noalign{\smallskip}
f_2(z)=z_2-b^2 z_1z_2-b z_2^2+z_1^3,\cr}
$$
for any $b\in\C^*$, and it is easy to build similar examples of any
order. At present, the theory in this case is satisfactorily developed for $n=2$
only. In particular, in [A2] is proved the following

\newthm Theorem \Abate: (Abate, 2001 [A2]) Every holomorphic local dynamical
system $f\in\End(\C^2,O)$ tangent to the identity, with an isolated fixed point,
admits at least one Fatou flower tangent to some direction.

\Rem Bracci and Suwa have proved a version of Theorem~\Abate\ for
$f\in\End(M,p)$ where $M$ is a {\it singular} variety with not too bad a
singularity at~$p$; see [BrS] for details. 

Let us describe the main ideas in the proof of Theorem~\Abate, because they
provide some insight on the dynamical structure of holomorphic local
dynamical systems tangent to the identity, and on how to deal with it. A shorter proof, 
deriving this theorem directly from Camacho-Sad theorem [CS] on the existence of
separatrices for holomorphic vector fields in $\C^2$,
is presented in [BCL] (see also [DI2]); however, such an approach provides fewer informations on 
the dynamical and geometrical structures of local dynamical systems
tangent to the identity.

The first idea is to exploit in a systematic way the transformation~\eqblowup,
following a procedure borrowed from algebraic geometry. 

\Def If $p$ is a point in a
complex manifold~$M$, there is a canonical way (see, e.g., [GH] or [A1]) to build a complex
manifold~$\tilde M$, called the {\sl blow-up} of~$M$ at~$p$, provided with a
holomorphic projection $\pi\colon\tilde M\to M$, so that $E=\pi^{-1}(p)$,
the {\sl exceptional divisor} of the blow-up, is canonically biholomorphic
to~$\P(T_pM)$, and $\pi|_{\tilde M\setminus E}\colon\tilde M\setminus E\to
M\setminus\{p\}$ is a biholomorphism. In suitable local coordinates, the
map~$\pi$ is exactly given by~\eqblowup. Furthermore, if $f\in\End(M,p)$ is
tangent to the identity, there is a unique way to lift~$f$ to a map~$\tilde
f\in\End(\tilde M,E)$ such that $\pi\circ\tilde f=f\circ\pi$, where $\End(\tilde
M,E)$ is the set of holomorphic maps defined in a neighbourhood of~$E$ with
values in~$\tilde M$ and which are the identity on~$E$. 

In particular, the
characteristic directions of~$f$ become points in the domain of the lifted map~$\tilde f$; 
and we shall see that this approach allows to determine which characteristic directions are
dynamically meaningful.

The blow-up procedure reduces the study of the dynamics of local holomorphic dynamical
systems tangent to the identity to the study of the dynamics of germs $f\in\End(M,E)$,
where $M$ is a complex $n$-dimensional manifold, and $E\subset M$ is a compact smooth
complex hypersurface pointwise fixed by~$f$. In [A2], [BrT] and [ABT1] we discovered
a rich geometrical structure associated to this situation.

Let $f\in\End(M,E)$ and 
take $p\in E$. Then for every
$h\in\ca O_{M,p}$ (where $\ca O_M$ is the structure sheaf of~$M$) the germ $h\circ f$ is
well-defined, and we have $h\circ f-h\in\ca I_{E,p}$, where $\ca I_E$ is the ideal sheaf 
of~$E$. 

\Def The {\sl $f$-order of vanishing} at~$p$ of~$h\in\ca O_{M,p}$ is
$$
\nu_f(h;p)=\max\{\mu\in\N\mid h\circ f-h\in\ca I_{E,p}^\mu\}\;,
$$
and the {\sl order of contact}~$\nu_f$ of~$f$ with~$E$ is
$$
\nu_f=\min\{\nu_f(h;p)\mid h\in\ca O_{M,p}\}\;.
$$

In [ABT1] we proved that $\nu_f$ does not depend on $p$, and that 
$$
\nu_f=\min_{j=1,\ldots,n} \nu_f(z^j;p)\;,
$$
where $(U,z)$ is any local chart centered at~$p\in E$ and $z=(z^1,\ldots,z^n)$. In particular, if the local chart $(U,z)$ is such that $E\cap U=\{z^1=0\}$ (and we shall say that the 
local chart is {\sl adapted} to $E$)
then setting $f^j=z^j\circ f$ we can write
$$
f^j(z)=z^j+(z^1)^{\nu_f}g^j(z)\;,
\neweq\eqquno
$$
where at least one among $g^1,\ldots, g^n$ does not vanish identically 
on $U\cap E$. 

\Def A map $f\in\End(M,E)$ is {\sl tangential} to~$E$ if 
$$
\min\bigl\{\nu_f(h;p)\mid h\in\ca I_{E,p}\bigr\}>\nu_f
$$
for some (and hence any) point $p\in E$. 

Choosing a local chart $(U,z)$ adapted to~$E$ so that we can express the coordinates of~$f$ in the form~\eqquno, it turns out that $f$ is tangential if and only if $g^1|_{U\cap E}\equiv 0$.

The $g^j$'s in \eqquno\ depend in general on the chosen chart; however, in [ABT1] we proved that setting
$$
\ca X_f= \sum_{j=1}^n g^j
{\de\over\de z^j}\otimes (dz^1)^{\otimes\nu_f}
\neweq\eqcaXf
$$
then $\ca X_f|_{U\cap E}$ defines a {\it global} section $X_f$ of the bundle $TM|_E\otimes (N_E^*)^{\otimes\nu_f}$, where $N_E^*$ is the conormal bundle 
of~$E$ into~$M$. 
The bundle $TM|_E\otimes (N_E^*)^{\otimes\nu_f}$ is canonically isomorphic to the bundle~$\Hom(N_E^{\otimes\nu_f},TM|_E)$. Therefore the section~$X_f$ induces a morphism still denoted by $X_f\colon N_E^{\otimes\nu_f}\to TM|_E$.

\Def The morphism $X_f\colon N_E^{\otimes\nu_f}\to TM|_E$ just defined
is the {\sl canonical morphism} associated to $f\in\End(M,E)$.

\Rem It is easy to check that $f$ is tangential if and only if the image of $X_f$ is contained in~$TE$. Furthermore, if $f$ is the lifting of a germ $f_o\in\End(
\C^n,O)$ tangent to the identity, then (see [ABT1]) $f$ is tangential if and only if $f_o$ is
non-dicritical; so in this case tangentiality is generic. Finally, in [A2] we used the term
``non degenerate" instead of "tangential".

\Def Assume that $f\in\End(M,E)$ is tangential.
We shall say that $p\in E$ is a {\sl singular point for~$f$} if $X_f$ vanishes at~$p$. 

By definition, $p\in E$ is a singular point for~$f$ if and only if
$$
g^1(p)=\cdots=g^n(p)=0
$$
for any local chart adapted to~$E$; so singular points are generically isolated.

In the tangential case, only
singular points are dynamically meaningful. Indeed, a not too difficult
computation (see~[A2], [AT1] and [ABT1])
yields the following

\newthm Proposition \sing: Let $f\in\End(M,E)$ be tangential,
and take $p\in E$. If $p$ is not singular,
then the stable set of the germ of $f$ at~$p$ coincides with $E$.

The notion of singular point allows us to identify the 
dynamically meaningful characteristic
directions. 

\Def Let $M$ be the blow-up of~$\C^n$ at the origin, and~$f$ the lift of 
a non-dicritical holomorphic local dynamical system $f_o\in\End(\C^n,O)$ 
tangent to the identity.
We shall say that $[v]\in\P^{n-1}(\C)=E$ is a {\sl singular direction} of~$f_o$ if it is a singular point of~$\tilde f$. 

It turns out that non-degenerate
characteristic directions are always singular (but the converse does not
necessarily hold), and that singular directions are always characteristic (but
the converse does not necessarily hold). Furthermore, the singular directions are the
dynamically interesting characteristic directions, because 
Propositions \suno\ and \sing\ imply that if $f_o$ has a non-trivial orbit
converging to the origin tangentially to~$[v]\in\P^{n-1}(\C)$ then $[v]$
must be a singular direction.

The important feature of the blow-up procedure is that, even though the underlying
manifold becomes more complex, the lifted maps become simpler. Indeed, using
an argument similar to one (described, for instance, in~[MM]) used in the study
of singular holomorphic foliations of 2-dimensional complex manifolds,
in~[A2] it is shown that after a finite number of blow-ups our original
holomorphic local dynamical system~$f\in\End(\C^2,O)$ tangent
to the identity can be lifted to a
map~$\tilde f$ whose singular points (are finitely many and) satisfy one of
the following conditions:
\smallskip
\item{(o)} they are dicritical; or,
\item{($\star$)}in suitable local coordinates centered at the singular point we
can write
$$
\cases{\tilde f_1(z)=z_1+\ell(z)\bigl(\lambda_1 z_1+O(\|z\|^2)\bigr),\cr
\tilde f_2(z)=z_2+\ell(z)\bigl(\lambda_2 z_2+O(\|z\|^2)\bigr),
\cr}
$$
with
\itemitem{($\star_1$)} $\lambda_1$,~$\lambda_2\ne 0$ and
$\lambda_1/\lambda_2$,~$\lambda_2/\lambda_1\notin\N$, or
\itemitem{($\star_2$)} $\lambda_1\ne 0$, $\lambda_2=0$.

\Rem This ``reduction of the singularities" statement holds only in dimension~2,
and it is not clear how to replace it in higher dimensions.

It is not too difficult to prove that Theorem~\EcalleHakim\ can be applied to both dicritical and
$(\star_1)$ singularities; therefore as soon as this blow-up procedure produces
such a singularity, we get a Fatou flower for the original dynamical system~$f$.

So to end the proof of Theorem~\Abate\ it suffices to prove that any such
blow-up procedure {\it must} produce at least one dicritical or $(\star_1)$
singularity. To get such a result, we need another ingredient.

Let again $f\in\End(M,E)$, where $E$ is a smooth complex hypersurface
in a complex manifold $M$, and assume that $f$ is tangential; let $E^o$
denote the complement in $E$ of the singular points of $f$. For simplicity of exposition we shall assume $\dim M=2$ and $\dim E=1$; but this part
of the argument works for any $n\ge 2$ (even when $E$ has singularities,
and it can also be adapted to non-tangential germs). 

Since $\dim E=1=\rk N_E$, the restriction of the canonical morphism $X_f$ to
$N_{E^o}^{\otimes\nu_f}$ is an isomorphism between $N_{E^o}^{\otimes\nu_f}$
and $TE^o$. Then
in [ABT1] we showed that it is possible to define a holomorphic 
connection~$\nabla$ on~$N_{E^o}$ by setting
$$
\nabla_u(s)=\pi([\ca X_f(\tilde u),\tilde s]|_S)\;,
\neweq\eqCSuno
$$
where: $s$ is a local section of $N_{E^o}$; $u\in TE^o$; $\pi\colon 
TM|_{E^o}\to N_{E^o}$ is the canonical projection; $\tilde s$ is any local
section of $TM|_{E^o}$ such
that $\pi(\tilde s|_{S^o})=s$; $\tilde u$ is any local section of
$TM^{\otimes\nu_f}$ such that $X_f\bigl(\pi(\tilde
u|_{E^o})\bigr)=u$; and $\ca X_f$ is locally given by \eqcaXf. In a chart 
$(U,z)$ adapted to~$E$, a local generator of $N_{E^o}$ is $\de_1=\pi(\de/\de z^1)$, a local generator of~$N_{E^o}^{\otimes\nu_f}$ is~$\de_1^{\otimes\nu_f}=\de_1\otimes\cdots\otimes\de_1$, and we have
$$
X_f(\de_1^{\otimes\nu_f})= g^2|_{U\cap E}{\de\over\de z^2}\;;
$$
therefore
$$
\nabla_{\de/\de z^2}\de_1=-\left.{1\over g^2}{\de g^1\over\de z^1}
\right|_{U\cap E}\de_1\;.
$$
In particular, $\nabla$ is a meromorphic connection on $N_E$, with poles
in the singular points of~$f$.

\Def The {\sl index} $\iota_p(f,E)$ of $f$ along $E$ at a point $p\in E$
is by definition the opposite of the residue at~$p$ of the connection~$\nabla$:
$$
\iota_p(f,E)=-\hbox{\rm Res}_p(\nabla)\;.
$$
In particular, $\iota_p(f,E)=0$ if $p$ is not a singular point of $f$.

\Rem If $[v]$ is a non-degenerate characteristic direction of a 
non-dicritical $f_o\in\End(\C^2,O)$ with non-zero director $\alpha\in\C^*$,
then it is not difficult to check that
$$
\iota_{[v]}(f,E)={1\over \alpha}\;,
$$
where $f$ is the lift of~$f_o$ to the blow-up of the origin.

Then in [A2] we proved the following {\sl index theorem} (see [Br1], 
[BrT] and [ABT1, 2]
for multidimensional versions and far reaching generalizations):

\newthm Theorem \index: ([A2], [ABT1]) Let $E$ be a compact Riemann surface
inside a 2-dimensional complex manifold~$M$. Take $f\in\End(M,E)$, and assume that~$f$ is
tangential to~$E$. Then
$$
\sum_{q\in E}\iota_q(f,E)=c_1(N_E)\;,
$$
where $c_1(N_E)$ is the first Chern class of the normal bundle~$N_E$ of~$E$
in~$M$.

Now, a combinatorial argument (inspired by Camacho and Sad~[CS]; see also [Ca] and [T]) shows
that if we have $f\in\End(\C^2,O)$ tangent to the identity with an isolated fixed point, and such that
applying the reduction of singularities to the lifted map~$\tilde f$ starting from a
singular direction~$[v]\in\P^1(\C)=E$ we end up only with $(\star_2)$
singularities, then the index of~$\tilde f$ at~$[v]$ along~$E$ must be a
non-negative rational number. But the first Chern class of~$N_E$ is~$-1$; so
there must be at least one singular directions whose index is not a non-negative
rational number. Therefore the reduction of singularities must yield at least one
dicritical or $(\star_1)$ singularity, and hence a Fatou flower for our 
map~$f$, completing the proof of Theorem~\Abate.

Actually, we have proved the following slightly more precise result:

\newthm Theorem \Abatedue: ([A2]) Let $E$ be a (not necessarily compact) Riemann
surface inside a 2-dimensional complex manifold~$M$, and take 
$f\in\End(M,E)$ tangential to~$E$. Let $p\in E$
be a singular point of~$f$ such that $\iota_p(f,E)\notin\Q^+$. Then there exist a Fatou flower for~$f$ at~$p$. In particular, if
$f_o\in\End(\C^2,O)$ is a non-dicritical
holomorphic local dynamical system tangent to the identity with an isolated
fixed point at the origin, and $[v]\in\P^1(\C)$ is a singular direction such
that~$\iota_{[v]}\bigl(f,\P^1(\C)\bigr)\notin\Q^+$, where $f$ is
the lift of~$f_o$ to the blow-up of the origin, then $f_o$ has a Fatou flower
tangent to~$[v]$.

\Rem This latter statement has been 
generalized in two ways. Degli Innocenti [DI1] has proved that we can allow~$E$
to be singular at~$p$ (but irreducible; in the reducible case one has to impose
conditions on the indeces of~$f$ along all irreducible components of~$E$ passing
through~$p$).  Molino [Mo], on the other hand, has proved that the statement
still holds assuming only $\iota_p(f,E)\ne 0$, at least for $f$ of order~2 (and
$E$ smooth at~$p$); it is natural to conjecture that this should be true for $f$
of any order.

As already remarked, the reduction of singularities via
blow-ups seem to work only in dimension~2. This leaves open the problem of the
validity of something like Theorem~\Abate\ in dimension~$n\ge 3$; see [AT1]
and [Ro2] for some partial results.

As far as I know, it is widely open, even in dimension~2, the
problem of describing the stable set of a holomorphic local dynamical
system tangent to the identity, as well as the more general problem of the
topological classification of such dynamical systems. Some results in
the case of a dicritical singularity are presented in~[BM]; for
non-dicritical singularities a promising 
approach in dimension~2 is described in~[AT3]. 

Let $f\in\End(M,E)$, where $E$ is a smooth Riemann surface
in a 2-dimensional complex manifold $M$, and assume that $f$ is tangential; let $E^o$
denote the complement in $E$ of the singular points of $f$. The connection~$\nabla$
on $N_{E^o}$ described above induces a connection (still denoted 
by $\nabla$) on $N_{E^o}^{\otimes \nu_f}$.

\Def In this setting, a {\sl geodesic} is a curve $\sigma\colon I\to E^o$
such that
$$
\nabla_{\sigma'}X_f^{-1}(\sigma')\equiv O\;.
$$

It turns out that $\sigma$ is a geodesic if and only if 
the curve $X_f^{-1}(\sigma')$ is an integral curves of a global
holomorphic vector field $G$ on the total space of $N_{E^o}^{\otimes\nu_f}$;
furthermore, $G$ extends holomorphically to the total space 
of $N_E^{\otimes\nu_f}$.

Now, assume that $M$ is the blow-up of the origin in~$\C^2$,
and $E$ is the exceptional divisor. Then there exists a canonical 
$\nu_f$-to-1 holomorphic covering map
$\chi_{\nu_f}\colon\C^2\setminus\{O\}\to N_E^{\otimes\nu_f}\setminus E$.
Moreover, if $f$ is the lift of a non-dicritical $f_o\in\End(\C^2,O)$ of
the form~\eqpsc\ with $P_\nu=(P_\nu^1,P_\nu^2)$, then $\nu_f=\nu-1$ and it turns out that $\chi_{\nu_f}$
maps integral curves of the homogeneous vector field 
$$
Q_\nu=P_\nu^1{\de\over\de z^1}+P_\nu^2{\de\over\de z^2}
$$
onto
integral curves of~$G$. In particular, to study the dynamics
of the time-1 map (which is tangent to the identity and of the form \eqpsc)of a non-dicritical homogeneous vector field~$Q_\nu$
it suffices to study the dynamics of such a geodesic vector field~$G$.
This is done in [AT3]; in particular, we get a complete description
of the local dynamics in a full neighbourhood of the
origin for a large class of holomorphic local dynamical systems tangent
to the identity. Since results like Theorem~\Camacho\ seems to
suggest that generic holomorphic local dynamical systems tangent
to the identity might be topologically conjugated to the time-1
map of a homogeneous vector field, this approach might eventually lead
to a complete topological description of the dynamics for generic
holomorphic local dynamical systems tangent
to the identity in dimension~2.

We end this section with a couple of words on holomorphic local dynamical
systems with a parabolic fixed point where the differential is not
diagonalizable. Particular examples are studied in detail in~[CD], [A4]
and~[GS]. In~[A1] it is described
a canonical procedure for lifting an~$f\in\End(\C^n,O)$ whose differential at the
origin is not diagonalizable to a map defined in a suitable iterated blow-up of
the origin (obtained blowing-up not only points but more general submanifolds)
with a canonical fixed point where the differential is diagonalizable. Using
this procedure it is for instance possible to prove the following

\newthm Corollary \Jordan: ([A2]) Let $f\in\End(\C^2,O)$ be a
holomorphic local dynamical system with $df_O=J_2$, the canonical Jordan matrix
associated to the eigenvalue~$1$, and assume that the origin is an isolated
fixed point. Then $f$ admits at least one parabolic curve tangent to~$[1:0]$ at
the origin. 

\smallsect 7. Several complex variables: other cases

Outside the hyperbolic and parabolic cases, there are not that many general
results. Theorems~\Bryunohyp\ and~\BryunoCremer\ apply to the elliptic case too,
but, as already remarked, it is not known whether the Bryuno condition is still
necessary for holomorphic linearizability. However,
another result in the spirit of Theorem~\BryunoCremer\ is the following:

\newthm Theorem \YC: (Yoccoz, 1995 [Y2]) Let $A\in GL(n,\C)$ be an invertible
matrix such that its eigenvalues have no resonances and such that its Jordan
normal form contains a non-trivial block associated to an eigenvalue of modulus
one. Then there exists $f\in\End(\C^n,O)$ with $df_O=A$ which is not
holomorphically linearizable.

A case that has received some attention is the so-called
semi-attractive case

\Def A holomorphic local dynamical system $f\in\End(\C^n,O)$
is said {\sl semi-attractive} if the eigenvalues of~$df_O$ are either equal to~1
or with modulus strictly less than~1. 

The dynamics of semi-attractive dynamical
systems has been studied by Fatou~[F4], Nishimura~[N],
Ueda~[U1--2], Hakim~[H1] and Rivi~[Ri--2]. Their results more or less
say that the eigenvalue~1 yields the existence of a ``parabolic manifold"~$M$ ---
in the sense of Theorem~\Hakim.(ii) --- of a suitable dimension, while the
eigenvalues with modulus less than one ensure, roughly speaking, that the
orbits of points in the normal bundle of~$M$ close enough to~$M$ are attracted
to it. For instance, Rivi proved the following

\newthm Theorem \Rivi: (Rivi, 1999 [Ri1--2]) Let $f\in\End(\C^n,O)$ be a
holomorphic local dynamical system. Assume that $1$ is an eigenvalue of
(algebraic and geometric) multiplicity~$q\ge 1$ of $df_O$, and that the other
eigenvalues of~$df_O$ have modulus less than~$1$. Then:
{\smallskip
\item{\rm(i)}We can choose local coordinates $(z,w)\in\C^q\times\C^{n-q}$ such
that $f$ expressed in these coordinates becomes
$$
\cases{f_1(z,w)=A(w)z+P_{2,w}(z)+P_{3,w}(z)+\cdots,\cr
\noalign{\smallskip}
f_2(z,w)=G(w)+B(z,w)z,\cr}
$$
where: $A(w)$ is a $q\times q$ matrix with entries holomorphic in~$w$ and
$A(O)=I_q$; the $P_{j,w}$ are $q$-uples of homogeneous polynomials in~$z$ of
degree~$j$ whose coefficients are holomorphic in~$w$; $G$ is a holomorphic
self-map of~$\C^{n-q}$ such that $G(O)=O$ and the eigenvalues of~$dG_O$ are
the eigenvalues of~$df_O$ with modulus strictly less than~$1$; and $B(z,w)$ is
an $(n-q)\times q$ matrix of holomorphic functions vanishing at the origin. In
particular, $f_1(z,O)$ is tangent to the identity.
\item{\rm(ii)} If $v\in\C^q\subset\C^m$ is a non-degenerate characteristic
direction for~$f_1(z,O)$, and the latter map has order~$\nu$,\break\indent then
there exist
$\nu-1$ disjoint $f$-invariant $(n-q+1)$-dimensional complex
submanifolds~$M_j$ of~$\C^n$, with\break\indent the origin in their boundary,
such that the orbit of every point of~$M_j$ converges to the origin tangentially
\break\indent to~$\C v\oplus E$, where $E\subset\C^n$ is the subspace generated
by the generalized eigenspaces associated to the\break\indent eigenvalues
of~$df_O$ with modulus less than one.}  

Rivi also has results in the spirit of Theorem~\Hakim, and results when the
algebraic and geometric multiplicities of the eigenvalue~$1$ differ; see~[Ri1, 2]
for details.

As far as I know, the only results on the formal or holomorphic
classification of semi-attractive holomorphic local dynamical systems
are due to Jenkins [J]. However,
building on work done by Canille Martins [CM] in dimension 2, and
using Theorem~\Camacho\ and general results on normally hyperbolic dynamical
systems due to Palis and Takens~[PT], Di Giuseppe has obtained the topological
classification when the eigenvalue $1$ has multiplicity 1 and the other
eigenvalues are not resonant:

\newthm Theorem \DiGiuseppe: (Di Giuseppe, 2004 [Di]) Let
$f\in\End(\C^n,O)$ be a holomorphic local dynamical system such that $df_O$ has
eigenvalues $\lambda_1$,~$\lambda_2,\ldots,\lambda_n\in\C$, where $\lambda_1$ is
a primitive $q$-root of unity, and $|\lambda_j|\ne 0$,~$1$ for $j=2,\ldots,n$.
Assume moreover that $\lambda_2^{r_2}\cdots\lambda_n^{r^n}\ne 1$ for all
multi-indeces $(r_2,\ldots, r_n)\in\N^{n-1}$ such that $r_2+\cdots+r_n\ge 2$.
Then $f$ is topologically locally conjugated either to $df_O$ or to the map
$$
z\mapsto (\lambda_1 z_1+z_1^{kq+1},\lambda_2 z_2,\ldots,\lambda_n z_n)
$$
for a suitable $k\in\N^*$.

We end this survey by recalling results by Bracci and Molino, and by Rong. Assume
that $f\in\End(\C^2,O)$ is a holomorphic local dynamical system such that the
eigenvalues of~$df_O$ are~$1$ and~$e^{2\pi i\theta}\ne 1$. If $e^{2\pi i\theta}$
satisfies the Bryuno condition, P\"oschel~[P\"o] proved the existence of a
1-dimensional $f$-invariant holomorphic disk containing the origin where~$f$ is
conjugated to the irrational rotation of angle~$\theta$. On the other hand,
Bracci and Molino give sufficient conditions (depending on~$f$ but not
on~$e^{2\pi i\theta}$, expressed in terms of two new holomorphic invariants, and
satisfied by generic maps) for the existence of parabolic curves tangent to the
eigenspace of the eigenvalue~1; see [BrM] for details, and [Ro3] for 
generalizations to $n\ge 3$.

\setref{EHRS}
\beginsection References

\art A1 M. Abate: Diagonalization of non-diagonalizable discrete holomorphic
dynamical systems! Amer. J. Math.! 122 2000 757-781

\art A2 M. Abate: The residual index and the dynamics of holomorphic maps
tangent to the identity! Duke Math. J.! 107 2001 173-207

\book A3 M. Abate: An introduction to hyperbolic dynamical systems! I.E.P.I.
Pisa, 2001

\art A4 M. Abate: Basins of attraction in quadratic dynamical systems with a
Jordan fixed point! Nonlinear Anal.! 51 2002 271-282

\art AT1 M. Abate, F. Tovena: Parabolic curves in $\C^3$! 
Abstr. Appl. Anal.! 2003 2003 275-294

\art AT2 M. Abate, F. Tovena: Formal classification of holomorphic maps tangent
to the identity! Disc. Cont. Dyn. Sys.! Suppl. 2005 1-10

\pre AT3 M. Abate, F. Tovena: Poincar\'e-Bendixson theorems for meromorphic connections
and homogeneous vector fields! Preprint! Institut Mittag-Leffler, 2008 

\art ABT1 M. Abate, F. Bracci, F. Tovena: Index theorems for holomorphic
self-maps! Ann. of Math.! 159 2004 819-864

\pre ABT2 M. Abate, F. Bracci, F. Tovena: Index theorems for holomorphic maps and
foliations! To appear in Indiana Univ. Math. J.! {\bf 57} (2008)

\art AM A. Abbondandolo, P. Majer: On the global stable manifold!  Studia Math.! 177  2006
113-131

\book Ar V.I. Arnold: Geometrical methods in the theory of ordinary differential
equations! Springer-Verlag, Berlin, 1988

\art B1 K. Biswas: Smooth combs inside hedgehogs!  Discrete Contin. Dyn. Syst.!  12  2005  853-880 

\pre B2 K. Biswas: Complete conjugacy invariants of nonlinearizable holomorphic dynamics!
Preprint, University of Pisa! 2007

\art B\"o L.E. B\"ottcher: The principal laws of convergence of iterates and their 
application to analysis! Izv. Kazan. Fiz.-Mat. Obshch.! 14 1904
155-234

\art Br1 F. Bracci: The dynamics of holomorphic maps near curves of fixed points!
Ann. Scuola Norm. Sup. Pisa! 2 2003 493-520

\art Br2 F. Bracci: Local dynamics of holomorphic diffeomorphisms! Boll. UMI ! 7--B 2004 609-636

\art BrM F. Bracci, L. Molino: The dynamics near quasi-parabolic fixed points of 
holomorphic diffeomorphisms in $\C^2$! Amer. J. Math.! 126 2004
671-686

\art BrS F. Bracci, T. Suwa: Residues for singular pairs and dynamics of
biholomorphic maps of singular surfaces! Internat. J. Math.! 15 2004
443-466 

\art BrT F. Bracci, F. Tovena: Residual indices of
holomorphic maps relative to  singular curves of fixed points on
surfaces! Math. Z.! 242 2002 481-490

\art BM F.E. Brochero Mart\'\i nez: Groups of germs of analytic diffeomorphisms
in $(\C^2,O)$! J. Dynamic. Control Systems! 9 2003 1-32

\art BCL F.E. Brochero Mart\'\i nez, F. Cano, L. L\'opez-Hernanz: Parabolic curves for diffeomorphisms in $\C^2$!  Publ. Mat.! 52  2008 189-194

\art Bry1 A.D. Bryuno: Convergence of transformations of differential equations
to normal forms! Dokl. Akad. Nauk. USSR! 165 1965 987-989

\art Bry2 A.D. Bryuno: Analytical form of differential equations, I! Trans.
Moscow  Math. Soc.! 25 1971 131-288

\art Bry3 A.D. Bryuno: Analytical form of differential equations, I\negthinspace
I! Trans. Moscow  Math. Soc.! 26 1972 199-239

\art C C. Camacho: On the local structure of conformal mappings and holomorphic
vector fields! Ast\'e\-risque! 59--60 1978 83-94

\book BH X. Buff, J.H. Hubbard: Dynamics in one complex variable!
To be published by Matrix Edition, Ithaca, NY

\art CS C. Camacho, P. Sad: Invariant varieties through singularities of holomorphic vector
fields! Ann. of Math.! 115 1982 579-595

\art CM J.C. Canille Martins: Holomorphic flows in $(\C^3,O)$ with resonances!
Trans. Am. Math. Soc.! 329 1992 825-837

\art Ca J. Cano: Construction of invariant curves for singular holomorphic
vector fields! Proc. Am. Math. Soc.! 125 1997 2649-2650

\book CG S. Carleson, F. Gamelin: Complex dynamics! Springer, Berlin, 1994

\book Ch M. Chaperon: G\'eom\'etrie diff\'erentielle et singularit\'es des
syst\`emes dynamiques! Ast\'e\-risque {\bf 138--139,} 1986

\art CD D. Coman, M. Dabija: On the dynamics of some diffeomorphisms of $\C^2$
near parabolic fixed points! Houston J. Math.! 24 1998 85-96

\art Cr1 H. Cremer: Zum Zentrumproblem! Math. An..! 98 1927 151-163

\art Cr2 H. Cremer: \"Uber die H\"aufigkeit der Nichtzentren! Math. Ann.! 115
1938 573-580

\art DI1 F. Degli Innocenti: Holomorphic dynamics near germs of singular curves!
 Math. Z. ! 251  2005 943-958

\pre DI2 F. Degli Innocenti: On the relations between discrete and continuous dynamics
in~$\C^2$! Ph.D. Thesis, Universit\`a di Pisa! 2007

\art DG D. DeLatte, T. Gramchev: Biholomorphic maps with linear parts having
Jordan blocks: linearization and resonance type phenomena! Math. Phys. El. J.! 8
2002 1-27

\art Di P. Di Giuseppe: Topological classification of holomorphic,
semi-hyperbolic germs in quasi-absence of resonances!  Math. Ann.! 334  2006 609-625

\art DH A. Douady, J.H. Hubbard: On the dynamics of polynomial-like mappings!
Ann. Sc. \'Ec. Norm. Sup.! 18 1985 287-343

\art D H. Dulac: Recherches sur les points singuliers des \'equationes 
diff\'erentielles! J. \'Ec. Polytechnique! IX 1904 1-125

\art \'E1 J. \'Ecalle: Th\'eorie it\'erative: introduction \`a la th\'eorie des invariants holomorphes! J. Math. Pures Appl.! 54 1975 183-258
 
\book \'E2 J. \'Ecalle: Les fonctions r\'esurgentes. Tome
I: Les alg\`ebres de fonctions r\'esurgentes! Publ. Math. Orsay {\bf 81-05,}
Universit\'e de Paris-Sud, Orsay, 1981 

\book \'E3 J. \'Ecalle: Les fonctions r\'esurgentes. Tome
I\negthinspace I: Les fonctions r\'esurgentes appliqu\'ees \`a l'it\'eration!
Publ. Math. Orsay {\bf 81-06,} Universit\'e de Paris-Sud, Orsay, 1981 

\book \'E4 J. \'Ecalle: Les fonctions r\'esurgentes. Tome I\negthinspace
I\negthinspace I: L'\'equation du pont et la classification analytique des
objects locaux! Publ. Math. Orsay {\bf 85-05,} Universit\'e de Paris-Sud, Orsay,
1985 

\coll \'E5 J. \'Ecalle: Iteration and analytic classification of local
diffeomorphisms of $\C^\nu$! Iteration theory and its functional equations!
Lect. Notes in Math. {\bf 1163,} Springer-Verlag, Berlin, 1985, pp. 41--48

\art \'EV J. \'Ecalle, B. Vallet: Correction and linearization of resonant 
vector fields and diffeomorphisms! Math. Z.! 229 1998 249-318

\art F1 P. Fatou: Sur les \'equations fonctionnelles, I! Bull. Soc. Math. France!
47 1919 161-271

\art F2 P. Fatou: Sur les \'equations fonctionnelles, I\negthinspace I! Bull.
Soc. Math. France! 48 1920 33-94

\art F3 P. Fatou: Sur les \'equations fonctionnelles, I\negthinspace
I\negthinspace I! Bull. Soc. Math. France! 48 1920 208-314

\art F4 P. Fatou: Substitutions analytiques et \'equations fonctionnelles \`a
deux variables! Ann. Sc. Ec. Norm. Sup.! 40 1924 67-142

\art Fa C. Favre: Classification of 2-dimensional contracting rigid germs and Kato surfaces. I!  J. Math. Pures Appl.!  79  2000  475-514

\book FJ1 C. Favre, M. Jonsson: The valuative tree! Lect. Notes in Math.
1853, Springer Verlag, Berlin, 2004

\art FJ2 C. Favre, M. Jonsson: Eigenvaluations!  Ann. Sci. \'Ecole Norm. Sup.! 40 2007 309-349 

\coll FHY A. Fathi, M. Herman, J.-C. Yoccoz: A proof of Pesin's stable manifold
theorem! Geometric Dynamics! Lect Notes in Math. 1007, Springer Verlag, Berlin,
1983, pp.~177-216

\art GS V. Gelfreich, D. Sauzin: Borel summation and splitting of separatrices
for the H\'enon map! Ann. Inst. Fourier Grenoble! 51 2001 513-567

\book GH P. Griffiths, J. Harris: Principles of algebraic geometry! Wyley, New York, 1978

\art G1 D.M. Grobman: Homeomorphism of systems of differential equations! Dokl.
Akad. Nauk. USSR! 128 1959 880-881

\art G2 D.M. Grobman: Topological classification of neighbourhoods of a
singularity in $n$-space! Math. Sbornik! 56 1962 77-94

\art H J.S. Hadamard: Sur l'it\'eration et les solutions asymptotyques des
\'equations diff\'erentielles! Bull. Soc. Math. France! 29 1901 224-228

\art Ha1 M. Hakim: Attracting domains for semi-attractive transformations
of~$\C^p$! Publ. Matem.! 38 1994 479-499

\art Ha2 M. Hakim: Analytic transformations of $(\C^p,0)$ tangent to the
identity! Duke Math. J.! 92 1998 403-428

\pre Ha3 M. Hakim: Transformations tangent to the identity. Stable pieces
of manifolds! Preprint! 1997

\art Har P. Hartman: A lemma in the theory of structural stability of
differential equations! Proc. Am. Math. Soc.! 11 1960 610-620

\book HK B. Hasselblatt, A. Katok: Introduction to the modern theory of
dynamical systems! Cambridge Univ. Press, Cambridge, 1995

\coll He M. Herman: Recent results and some open questions on Siegel's
linearization theorem of germs of complex analytic diffeomorphisms of $\C^n$
near a fixed point! Proc. $8^{th}$ Int. Cong. Math. Phys.! World Scientific,
Singapore, 1986, pp. 138--198

\book HPS M. Hirsch, C.C. Pugh, M. Shub: Invariant manifolds! Lect. Notes Math.
{\bf 583,} Springer-Verlag, Berlin, 1977

\art HP J.H. Hubbard, P. Papadopol: Superattractive fixed points in $\C^n$!
Indiana Univ. Math. J.! 43 1994 321-365

\art I1 Yu.S. Il'yashenko: Divergence of series reducing an analytic differential
equation to linear normal form at a singular point! Funct. Anal. Appl.! 13 1979
227-229

\coll I2 Yu.S. Il'yashenko: Nonlinear Stokes phenomena! Nonlinear Stokes
phenomena! Adv. in Soviet Math. {\bf 14,} Am. Math. Soc., Providence, 1993,
pp.~1--55

\book IY Y.S. Il'yashenko, S. Yakovenko: Lectures on analytic differential
equations! Graduate Studies in Mathematics 86, American Mathematical
Society, Providence, RI, 2008

\art J A. Jenkins: Further reductions of Poincar\'e-Dulac normal forms
in $\C^{n+1}$! Proc. Am. Math. Soc.! 136 2008 1671-1680

\art K T. Kimura: On the iteration of analytic functions! Funk. Eqvacioj! 14
1971 197-238

\art {K\oe}  G. K\oe nigs: Recherches sur les integrals de certain equations
fonctionelles! Ann. Sci. \'Ec. Norm. Sup.! 1 1884 1-41

\art L L. Leau: \'Etude sur les equations fonctionelles \`a une ou plusieurs variables! Ann.
Fac. Sci. Toulouse! 11 1897 E1-E110

\art M1 B. Malgrange: Travaux d'\'Ecalle et de Martinet-Ramis sur les syst\`emes 
dynamiques! Ast\'erisque! 92-93 1981/82 59-73

\art M2 B. Malgrange: Introduction aux travaux de J. \'Ecalle! Ens. Math.!
31 1985 261-282

\book Ma S. Marmi: An introduction to small divisors problems! I.E.P.I., Pisa, 
2000

\art MM J.F. Mattei, R. Moussu: Holonomie et int\'egrales premi\`eres! Ann. Scient. Ec.
Norm. Sup.! 13 1980 469-523

\book Mi J. Milnor: Dynamics in one complex variable. Third edition! Annals of Mathematics Studies, 160, Princeton University Press, Princeton, NJ, 2006

\art Mo L. Molino: The dynamics of maps tangent to the identity and with nonvanishing
index! Trans. Amer. Math. Soc.! 361 2009 1597-1623

\art N Y. Nishimura: Automorphismes analytiques admettant des sousvari\'et\'es
de point fixes attractives dans la direction transversale! J. Math. Kyoto Univ.!
23 1983 289-299

\art PT J. Palis, F. Takens: Topological equivalence of normally hyperbolic
dynamical systems! Topology! 16 1977 335-345

\art P1 R. P\'erez-Marco: Sur les dynamiques holomorphes non lin\'earisables et
une conjecture de V.I. Arnold! Ann. Sci. \'Ecole Norm. Sup.! 26 1993 565-644

\pre P2 R. P\'erez-Marco: Topology of Julia sets and hedgehogs! Preprint!
Universit\'e de Paris-Sud, 94-48, 1994

\art P3 R. P\'erez-Marco: Non-linearizable holomorphic dynamics having an
uncountable number of symmetries! Invent. Math.! 199 1995 67-127

\pre P4 R. P\'erez-Marco: Holomorphic germs of quadratic type! Preprint! 1995

\pre P5 R. P\'erez-Marco: Hedgehogs dynamics! Preprint! 1995

\art P6 R. P\'erez-Marco: Sur une question de Dulac et Fatou! C.R. Acad. Sci.
Paris! 321 1995 1045-1048

\art P7 R. P\'erez-Marco: Fixed points and circle maps! Acta Math.! 179 1997
243-294

\pre P8 R. P\'erez-Marco: Linearization of holomorphic germs with resonant
linear part! Preprint, arXiv: math.DS/0009030! 2000

\art P9 R. P\'erez-Marco: Total convergence or general divergence in small
divisors! Comm. Math. Phys.! 223 2001 451-464 

\art Pe O. Perron: \"Uber Stabilit\"at und asymptotisches Verhalten der Integrale
von Differentialgleichungssystemen! Math. Z.! 29 1928 129-160

\art Pes Ja.B. Pesin: Families of invariant manifolds corresponding to non-zero
characteristic exponents! Math. USSR Izv.! 10 1976 1261-1305

\book Po H. Poincar\'e: \OE uvres, Tome I! Gauthier-Villars, Paris, 1928,
pp.~XXXVI--CXXIX

\art P\"o J. P\"oschel: On invariant manifolds of complex analytic mappings near
fixed points! Exp. Math.! 4 1986 97-109

\pre R J. Raissy: Linearization of holomorphic germs with quasi-Brjuno fixed points! Preprint, arXiv: math.DS/0710.3650v3! 2008

\art Re1 L. Reich: Das Typenproblem bei formal-biholomorphien Abbildungen mit
anziehendem Fixpunkt! Math. Ann.! 179 1969 227-250

\art Re2 L. Reich: Normalformen biholomorpher Abbildungen mit anziehendem
Fixpunkt! Math. Ann.! 180 1969 233-255

\book Ri1 M. Rivi: Local behaviour of discrete dynamical systems! Ph.D. Thesis,
Universit\`a di Firenze, 1999

\art Ri2 M. Rivi: Parabolic manifolds for semi-attractive holomorphic germs!
Mich. Math. J.! 49 2001 211-241

\art Ro1 F. Rong: Linearization of holomorphic germs with quasi-parabolic fixed points!
Ergodic Theory Dynam. Systems!  28  2008  979-986

\pre Ro2 F. Rong: Robust parabolic curves in $\C^m$ ($m\ge 3$)! To appear in Houston 
J. Math.! 2008 

\art Ro3 F. Rong: Quasi-parabolic analytic transformations of $\C^n$!  J. Math. Anal. Appl. !
343  2008  99-109

\art R\"u H. R\"ussmann: Stability of elliptic fixed points of analytic area-preserving mappings under the Brjuno condition! Ergodic Theory Dynam. Systems! 22 2002 1551-1573

\art S A.A. Shcherbakov: Topological classification of germs of conformal
mappings with identity linear part! Moscow Univ. Math. Bull.! 37 1982 60-65

\book Sh M. Shub: Global stability of dynamical systems! Springer, Berlin, 1987

\art Si C.L. Siegel: Iteration of analytic functions! Ann. of Math.! 43 1942
607-612

\art St1 S. Sternberg: Local contractions and a theorem of Poincar\'e! Amer. J.
Math.! 79 1957 809-824

\art St2 S. Sternberg: The structure of local homomorphisms, I\negthinspace I!
Amer. J. Math.! 80 1958 623-631

\art T M. Toma: A short proof of a theorem of Camacho and Sad!  Enseign. Math.!  45  1999  311-316

\art U1 T. Ueda: Local structure of analytic transformations of two complex
variables, I! J. Math. Kyoto Univ.! 26 1986 233-261

\art U2 T. Ueda: Local structure of analytic transformations of two complex
variables, I\negthinspace I! J. Math. Kyoto Univ.! 31 1991 695-711

\art Us S. Ushiki: Parabolic fixed points of two-dimensional complex dynamical
systems! S\=urikaisekiken\-ky\=usho K\=oky\=uroku! 959 1996 168-180

\art V S.M. Voronin: Analytic classification of germs of conformal maps
$(\C,0)\to (\C,0)$ with identity linear part! Func. Anal. Appl.! 15 1981 1-17

\art Y1  J.-C. Yoccoz: Lin\'earisation des germes de diff\'eomorphismes
holomorphes de $(\C,0)$! C.R. Acad. Sci. Paris! 306 1988 55-58

\art Y2 J.-C. Yoccoz: Th\'eor\`eme de Siegel, nombres de Bryuno et polyn\^omes
quadratiques! Ast\'e\-ris\-que! 231 1995 3-88

\art W B.J. Weickert: Attracting basins for automorphisms of $\C^2$! Invent.
Math.! 132 1998 581-605

\art Wu H. Wu: Complex stable manifolds of holomorphic diffeomorphisms! Indiana
Univ. Math. J.! 42 1993 1349-1358

\bye